\documentclass[11pt]{amsart}
\usepackage{latexsym}
\usepackage{amssymb}
\usepackage{pstricks,pst-node}

\psset{unit=1pt} \psset{arrowsize=4pt 1} \psset{linewidth=1pt}
\newrgbcolor{mygreen}{.2 .7 .2}
\newrgbcolor{myred}{.7 .3 .2}
\newgray{mygray}{.90}
\countdef\x=23 \countdef\y=24 \countdef\z=25 \countdef\t=26

\def\abox(#1,#2)#3{\psset{linewidth=1pt}
\x=#1 \y=#2 \multiply\x by 16 \multiply\y by 16 \z=\x \t=\y
\advance\z by 16 \advance\t by 16
\psline(\x,\y)(\x,\t)(\z,\t)(\z,\y)(\x,\y) \advance\x by 8
\advance\y by 8 \rput(\x,\y){{\bf #3}}}

\def\bbox(#1,#2)#3{ \psset{linewidth=0.5pt}
\x=#1 \y=#2 \multiply\x by 16 \multiply\y by 16 \z=\x \t=\y
\advance\z by 16 \advance\t by 16
\psline(\x,\y)(\x,\t)(\z,\t)(\z,\y)(\x,\y) \advance\x by 8
\advance\y by 8 \rput(\x,\y){{#3}}}

\textwidth=\textwidth \textheight=\textheight \headheight=9pt
\theoremstyle{plain}
\newtheorem{thm}{Theorem}
\newtheorem{lem}[thm]{Lemma}
\newtheorem{prop}[thm]{Proposition}
\newtheorem{lemma}[thm]{Lemma}
\newtheorem{corollary}[thm]{Corollary}

\newtheorem{conjecture}[thm]{Conjecture}

\theoremstyle{definition}
\newtheorem{definition}[thm]{Definition}

\newtheorem{remark}[thm]{Remark}

\newtheorem{problem}{Problem}
\newtheorem{example}{Example}

\newcommand{\set}[1]{ \left\{ #1 \right\} }
\newcommand{\ip}[1]{\left<#1\right>}
\def\c{{\mathcal C}}
\def\u{{\mathbf u}}
\def\U{{\mathcal U}}
\def\tS{\tilde{S}}
\def\tF{\tilde{F}}
\def\h{\tilde{h}}
\def\L{\Lambda}
\def\Ln{\Lambda^{(n)}}
\def\nL{\Lambda_{(n)}}
\def\Sn{S_{n^\lozenge}}
\def\Snn{S^{n^\lozenge}}

\def\a{\alpha}
\def\C{\mathbb C}
\def\ll{\lambda}
\def\p{\mathcal P}
\def\Z{\mathbb Z}
\def\P{\mathbb P}
\def\N{\mathbb N}
\def\ln{\mathrm {ln}}
\def\rn{\mathrm {rn}}
\def\pa{\mathrm{Par}}

\begin{document}
\title{Affine Stanley symmetric functions}
\author{Thomas Lam}
\email{thomasl@math.mit.edu}
\address{Department of Mathematics,
         M.I.T., Cambridge, MA 02139}
\date{January, 2005}
\begin{abstract}We define a new family $\tilde{F}_w(X)$ of
generating functions for $w \in \tilde{S}_n$ which are affine
analogues of Stanley symmetric functions.  We establish basic
properties of these functions including symmetry, dominance and
conjugation.  We conjecture certain positivity properties in terms
of a subfamily of symmetric functions called affine Schur
functions.  As applications, we show how affine Stanley symmetric
functions generalise the (dual of the) $k$-Schur functions of
Lapointe, Lascoux and Morse as well as the cylindric Schur
functions of Postnikov.  Conjecturally, affine Stanley symmetric
functions should be related to the cohomology of the affine flag
variety.
\end{abstract}

 \maketitle
 \section{Introduction}
In \cite{Sta84}, Stanley introduced a family $\set{F_w(X)}$ of
symmetric functions now known as \emph{Stanley symmetric
functions}.  He used these functions to study the number of
reduced decompositions of permutations $w \in S_n$.  Later, the
functions $F_w(X)$ were found to be closely related to the
Schubert polynomials of Lascoux and Sch\"{u}tzenberger
\cite{LS82}, which are well known to be related to the geometry of
flag varieties.

The aim of this paper is to define and study an analogue
$\tF_w(X)$ of Stanley symmetric functions for the affine symmetric
group $\tilde{S}_n$ which we call \emph{affine Stanley symmetric
functions}.  Our definition of $\tF_w(X)$ is motivated by
\cite{FS} and \cite{FG} and involve an algebra which we call the
\emph{affine nilCoxeter algebra}.  This algebra is an affine
version of the nilCoxeter algebra used in \cite{FS}.  When $w \in
S_n \subset \tS_n$, we have $\tF_w(X) = F_w(X)$. Our first main
theorem is that these functions $\tilde{F}_w(X)$ are indeed
symmetric functions. Imitating \cite{Sta84}, we show basic
properties of these functions:
\begin{enumerate}
\item
the relation to reduced words: \[ [x_1x_2\cdots x_{l(w)}]\tF_w(X)
= \#\{\mbox{reduced words of $w$}\},\]
\item
a skewing formula:
\[
s_1^\perp \cdot \tF_w = \sum_{w \gtrdot v} \tF_v
\]
where $\gtrdot$ denotes the covering relation in weak Bruhat
order,
\item
a conjugacy formula: \[ \tF_{w^*} = \omega^+(\tF_w) \] where
$*:\tS_n \rightarrow \tS_n$ and $\omega^+:\Ln \rightarrow \Ln$ are
involutions ($\Ln = \C\langle m_\ll \mid \; \ll_1 \leq
n-1\rangle$),
\item
and the existence of a unique dominant monomial term $m_{\mu(w)}$:
\[
\tF_w = m_{\mu(w)} + \sum_{\ll \prec \mu(w)} b_{w\ll} m_\ll.
\]
\end{enumerate}

An important special case occurs when $w$ is a \emph{Grassmannian}
permutation.  A permutation $w \in \tS_n$ is Grassmannian if it is
a minimal length coset representative of a coset of $S_n
\backslash \tS_n$.  In this case we obtain the \emph{affine Schur
functions} $\tF_\ll(X) = \tF_w(X)$ which may be labelled by
partitions with no part greater than $n-1$.  We show that the
affine Schur functions $\{\tF_\ll \mid \; \ll_1 \leq n-1\}$ form a
basis of the space $\Ln$ spanned by $\{m_\ll \mid \; \ll_1 \leq
n-1\}$ where the $m_\ll$ are monomial symmetric functions. Edelman
and Greene \cite{EG} and separately Lascoux and Sch\"{u}tzenberger
\cite{LS} have shown that Stanley symmetric functions $F_w(X)$
expand positively in terms of Schur functions $s_\ll(X)$.  We
conjecture that affine Stanley symmetric functions expand
positively in terms of affine Schur functions. We prove that a
unique maximal and minimal ``dominant'' term exists in such an
expansion.

Our definition of affine Stanley symmetric functions is motivated
by relations with two other classes of symmetric functions which
have received attention lately.  Lapointe, Lascoux and Morse
\cite{LLM} initiated the study of \emph{$k$-Schur functions},
denoted $s_\ll^{(k)}(X)$, in their study of Macdonald polynomial
positivity.  It is conjectured that $k$-Schur functions form an
``intermediate'' basis between the Macdonald polynomials
$\set{H_\mu(X;q,t)}$ \cite{Mac} and the Schur functions
$\set{s_\ll(X)}$ so that the transition coefficients are positive
in both intermediate steps. Lapointe and Morse have more recently
connected the multiplication of $k$-Schur functions with the
Verlinde algebra of $U(m)$.

Our affine Schur functions had earlier been defined using
``\emph{$k$-tableaux}'' by Lapointe and Morse, who called these
functions \emph{dual $k$-Schur functions}.  Work of Lapointe and
Morse \cite{LM04} relating $k$-Schur functions to $n$-cores and to
the affine symmetric group show that our definition of affine
Schur functions are indeed dual to $k$-Schur functions.  In this
context the symmetry of affine Schur functions is not obvious, but
follows from the symmetry of general affine Stanley symmetric
functions. The relation with $k$-Schur functions also suggest the
study of \emph{skew affine Schur functions} $\tF_{\ll/\mu}(X)$,
another special case of affine Stanley symmetric functions which
we study.

Separately, \emph{cylindric Schur functions} were defined by
Postnikov \cite{Pos} (see also \cite{GK}).  He showed that certain
coefficients of the expansion of \emph{toric Schur functions} (a
special case of cylindric Schur functions) in terms of Schur
functions were equal to the 3-point genus 0 \emph{Gromov-Witten
invariants} of the Grassmannian $Gr_{m,n}$ (which are the
multiplication constants of the \emph{quantum cohomology}
$QH^*(Gr_{m,n})$ of the Grassmannian). Cylindric Schur functions
are defined as generating functions of cylindric semi-standard
tableaux, which are tableaux drawn on a cylinder. We show that
cylindric Schur functions are special cases of skew affine Schur
functions and that they are exactly equal to affine Stanley
symmetric functions $\tF_{w}$ labelled by affine permutations $w$
which are ``$321$-avoiding''. These results are affine analogues
of some of the results in \cite{BJS}. However, the affine case is
significantly more difficult.  For example, any normal Stanley
symmetric function $F_w$ is equal to some skew affine Schur
function.

We also show that our conjecture that affine Stanley symmetric
functions expand positively in terms of affine Schur functions
implies both the Schur positivity of toric Schur polynomials and
the positivity of the multiplication of $k$-Schur functions (the
latter is currently still a conjecture).  Our work also explains
why cylindric Schur functions are not in general Schur positive in
infinitely many variables (see \cite{McN}).  The three families of
symmetric functions: affine Stanley, skew affine Schur and
cylindric Schur can be thought of as arising from three different
representations of the affine nilCoxeter algebra $\U_n$ in a
uniform manner.  The representation from which affine Stanley
symmetric functions arise is the left regular representation of
$\U_n$ and so leads to the most general symmetric functions which
can arise in this manner.

The connections with $k$-Schur functions and cylindric Schur
functions already indicate that affine Stanley symmetric functions
are important objects.  In fact our work, combined with the known
connection between quantum cohomology and the Verlinde algebra,
essentially implies that the main results of \cite{Pos} and
\cite{LM05} are equivalent.  However, it seems the most exciting
direction to take is to extend our definition of affine Stanley
symmetric functions $\tF_w$ to \emph{affine Schubert polynomials}
$\tilde{\mathfrak S}_w$ and connect them with the affine flag
variety ${\mathcal G}/{\mathcal B}$ of type $\tilde{A}_{n-1}$.
Usual Stanley symmetric functions are certain stable ``limits'' of
the Schubert polynomials ${\mathfrak S}_w$, which are well known
to represent Schubert varieties in the cohomology $H^*(G/B)$ of
the flag variety and possess numerous remarkable properties.  In
the affine case, the cohomology classes $[\Omega_w] \in
H^*({\mathcal G}/{\mathcal B})$ representing Schubert varieties
are labelled by $w \in \tS_n$ and should conjecturally be related
to affine Stanley symmetric functions.

In the Grassmannian case, Morse and Shimozono \cite{MS} have
conjectured that affine Schur functions represent the Schubert
classes in the cohomology of the affine Grassmannian. The study of
affine Schur functions should make explicit the relationship
between the affine Grassmannian $\mathcal{G}/\mathcal{P}$, and the
Verlinde algebras of $U(m)$ with level $n-m$ or the quantum
cohomology $QH^*(Gr_{m,n})$ of the Grassmannian, which are already
known to be connected; see for example \cite{Wit}.

Finally, we make the natural generalisation to \emph{affine stable
Grothendieck polynomials} $\tilde{G}_w(X)$ which speculatively
should be stable limits of the $K$-theory Schubert classes of the
affine flag variety.

Much work has also been done with a version of the Stanley
symmetric functions for the hyperoctahedral group; see
\cite{Kra,LTK,FK2}.  We intend to generalise this to the affine
case and also investigate Stanley symmetric functions for general
Coxeter groups in later work.

\bigskip

{\bf Overview.}  In sections \ref{sec:affineSymmetric} and
\ref{sec:symmetricFunctions}, we establish some notation for
affine permutations and for symmetric functions. In section
\ref{sec:stan} we recall the definition of Stanley symmetric
functions, give their main properties and explain the relationship
with Schubert polynomials. In section \ref{sec:affineStan}, we
define the affine nilCoxeter algebra and affine Stanley symmetric
functions and prove that the latter are symmetric.  In Section
\ref{sec:Vrep}, we explain how affine Stanley symmetric functions
arise from different representations of the affine nilCoxeter
algebra.  In section \ref{sec:basic}, we prove a coproduct formula
for affine Stanley symmetric functions.  In section
\ref{sec:monomialDominance}, we show that affine Stanley symmetric
functions have a unique dominant monomial term.  In section
\ref{sec:conjugacy}, we prove a conjugacy formula, imitating
\cite{Sta84}.  In section \ref{sec:affineSchur}, we define and
study affine Schur functions. In section \ref{sec:ncores}, we
study the relationship between $n$-cores and the affine symmetric
group, following in part \cite{Las}.  In sections \ref{sec:skew},
\ref{sec:ktableaux} and \ref{sec:kschur} we define skew affine
Schur functions and relate them to $k$-Schur functions. In section
\ref{sec:cylindric}, we recall the definition of a cylindric Schur
function and connect them with skew affine Schur functions.  In
section \ref{sec:321}, we show that cylindric Schur functions
correspond exactly to 321-avoiding permutations.  In section
\ref{sec:positivity}, we make a number of positivity conjectures
concerning the expansion of affine Stanley symmetric functions in
terms of affine Schur functions. Finally, in section
\ref{sec:final}, we discuss some further extensions of our theory
and in particular a generalisation to affine stable Grothendieck
polynomials.
\medskip

A condensed preliminary version of this paper appeared as
\cite{Lam1}.

\medskip

{\bf Acknowledgements.}  I thank Jennifer Morse for interesting
discussions about $k$-Schur functions and dual $k$-Schur
functions.  I also thank Mark Shimozono for discussions relating
to the affine Grassmannian. I am grateful to Alex Postnikov for
introducing cylindric and toric Schur functions in his class.  I
am also indebted to my advisor, Richard Stanley for guidance over
the last couple of years.

\section{Affine symmetric group}
\label{sec:affineSymmetric}

A positive integer $n \geq 3$ will be fixed throughout the paper.
Let $\tS_n$ denote the affine symmetric group with simple
generators $s_0,s_1,\ldots,s_{n-1}$ satisfying the relations
\begin{align*}
s_i s_{i+1} s_i &= s_{i+1} s_i s_{i+1} &  \mbox{for all $i$} \\
s_i^2 &= 1 &\mbox{for all $i$} \\  s_is_j &= s_j s_i &\mbox{for
$|i-j| \geq 2$}.
\end{align*}
Here and elsewhere, the indices will be taken modulo $n$ without
further mention.  One may realize $\tS_n$ as the set of all
bijections $w:\Z\rightarrow\Z$ such that $w(i+n)=w(i)+n$ for all
$i$ and $\sum_{i=1}^n w(i) = \sum_{i=1}^n i$. In this realization,
to specify an element $w \in \tS_n$ it suffices to give the
``window" $[w(1),w(2),\dotsc,w(n)]$.  The product $w\cdot v$ of
two affine permutations is then the composed bijection $w \circ v:
\Z \rightarrow \Z$.  Thus $ws_i$ is obtained from $w$ by swapping
the values of $w(i+kn)$ and $w(i+kn+1)$ for every $k \in \Z$.  See
\cite{BB} for more details.

The symmetric group $S_n$ embeds in $\tS_n$ as the subgroup
generated by $s_1,s_2, \ldots, s_{n-1}$.  Since there are many
embeddings of the $S_n$ into $\tS_n$ we will denote this
particular embedding by $\Sn$.

For an element $w \in \tS_n$ let $R(w)$ denote the set of
\emph{reduced words} for $w$.  A word $\rho =
(\rho_1\rho_2\cdots\rho_l) \in [0,n-1]^l$ is a reduced word for
$w$ if $w = s_{\rho_1} s_{\rho_2} \cdots s_{\rho_l}$ and $l$ is
the smallest possible integer for which such a decomposition
exists.  Abusing notation slightly, we also call $ s_{\rho_1}
s_{\rho_2} \cdots s_{\rho_l}$ a reduced word for $w$.
 The integer $l = l(w)$ is called the \emph{length} of $w$. If
$\rho, \pi \in R(w)$ for some $w$, then we write $\rho \sim \pi$.
If $\rho$ is an arbitrary word with letters from $[0,n-1]$ then we
write $\rho \sim 0$ if it is not a reduced word of any affine
permutation.  If $w, u \in \tS_n$ then we say that $w$
\emph{covers} $u$ if $w = s_i\cdot u$ and $l(w) = l(u) + 1$; and
we write $w \gtrdot u$. The transitive closure of $\gtrdot$ is
called the \emph{weak Bruhat order} and denoted $\geqslant$.

The \emph{code} $c(w)$ or \emph{affine inversion table} \cite{BB,
Las} of an affine permutation $w$ is a vector $c(w) =
(c_1,c_2,\ldots,c_{n}) \in \N^n - \P^n$ of non-negative entries
with at least one 0.  The entries are given by $c_i = \#\set{j \in
\Z \, \mid \, j > i \,\, \text{and} \,\, w(j) < w(i)}$.  It is
shown in \cite{BB} that there is a bijection between codes and
affine permutations and that $l(w) = |c(w)| = \sum_{i=1}^n c_i$.
The right action of the simple generator $s_i$ on the code $c =
(c_1,c_2,\ldots,c_{n})$ is given by
\[
(c_1,\ldots,c_{i},c_{i+1},\ldots,c_{n})\cdot s_i =
(c_1,\ldots,c_{i+1}+1,c_{i},\ldots,c_{n})
\]
whenever $c_i > c_{i-1}$.  Thus $c(ws_i) = c(w)\cdot s_i$.

\section{Symmetric functions}
\label{sec:symmetricFunctions} A partition $\ll = (\ll_1 \geq
\ll_2 \geq \cdots \geq \ll_l > 0)$ is a weakly decreasing finite
sequence of positive integers.  We use $\ll$, $\mu$ and $\nu$ to
denote partitions and will always draw them in the English
notation (top-left justified).  The \emph{dominance order}
$\preceq$ on partitions is given by $\ll \preceq \mu$ if and only
if $\ll_1 + \ll_2 + \cdots + \ll_i \leq \mu_1 +\mu_2 +\cdots +
\mu_i$ for every $i$.  We will also assume the reader is
reasonably familiar with the usual notions of corners, conjugates
and (semistandard) Young tableaux.

We will follow mostly \cite{Mac,EC2} for our symmetric function
notation.  Let $\Lambda$ denote the ring of symmetric functions
over $\C$. Usually, our symmetric functions will have an infinite
set of variables $x_1,x_2, \ldots$ and will be written as
$f(x_1,x_2,\ldots)$ or $f(X)$. If we need to emphasize the
variables used, we write $\Lambda_X$.

We will use $m_\ll$, $p_\ll$, $e_\ll$, $h_\ll$ and $s_\ll$ to
denote the \emph{monomial}, \emph{power sum}, \emph{elementary},
\emph{homogeneous} and \emph{Schur} bases of $\Lambda$.  It is
well known that $\{e_n\}$ and $\{h_n\}$ are algebraically
independent generators of $\Lambda$.  Let $\ip{.,.}$ denote the
\emph{Hall inner product} of $\Lambda$ satisfying
$\ip{h_\ll,m_\mu} = \ip{s_\ll,s_\mu} =\delta_{\ll\mu}$. For $f \in
\Lambda$, write $f^\perp: \Lambda \rightarrow \Lambda$ for the
linear operator adjoint to multiplication by $f$ with respect to
$\ip{.,.}$.  We let $\omega: \Lambda \rightarrow \Lambda$ denote
the $\C$-algebra involution of $\Lambda$ sending $h_n$ to $e_n$.

If $f(X) \in \Lambda$ then $f(x_1,x_2,\ldots,y_1,y_2,\ldots) =
f(X,Y) = \sum_i f_i(X) \otimes g_i(Y) \in \Lambda_X \otimes
\Lambda_Y$ for some $f_i$ and $g_i$. This is the \emph{coproduct}
of $f$, written $\Delta f = \sum_i f_i \otimes g_i \in \Lambda
\otimes \Lambda$. We have the following formula for the coproduct
(\cite{Mac}):
\begin{equation}
\label{eq:coproduct} \Delta f = \sum_\ll s_\ll^ \perp f \otimes
s_\ll. \end{equation}

The ring of symmetric functions $\Lambda$ is a self dual
Hopf-algebra with respect to $\ip{.,.}$, so that
\begin{equation}
\label{eq:hopf} \ip{\Delta f , g \otimes h} = \ip{f, gh},
\end{equation}
where $\ip{f_1 \otimes f_2, g_1\otimes g_2} :=
\ip{f_1,g_1}\ip{f_2,g_2}$.

Let $\pa^n$ denote the set $\set{\ll \mid \ll_1 \leq n-1}$ of
partitions with no row longer than $n-1$.  The following two
subspaces of $\Lambda$ will be important to us:
\begin{align*}
\Ln &:= \C \left< m_\ll \mid \ll \in \pa^n \right > \\
\nL &:= \C \left <h_\ll \mid \ll \in \pa^n \right> = \C \left
<e_\ll \mid \ll \in \pa^n \right> = \C \left <p_\ll \mid \ll \in
\pa^n \right>.
\end{align*}
If $f \in \nL$ and $g \in \Ln$ then define $\ip{f,g}$ to be their
usual Hall inner product within $\L$.  Thus $\set{h_\ll}$ and
$\set{m_\ll}$ with $\ll \in \pa^n$ form dual bases of $\nL$ and
$\Ln$.  Note that $\nL$ is a subalgebra of $\L$ but $\Ln$ is not
closed under multiplication.  Instead, $\Ln$ is a coalgebra; it is
closed under comultiplication.

\section{Stanley symmetric functions}
\label{sec:stan} Let $w \in S_n$ with length $l = l(w)$.  Define
the generating function $F_{w^{-1}}(X)$ by
\[
F_{w^{-1}}(x_1,x_2,\ldots) = \sum_{a_1a_2\cdots a_l \in R(w)}
\sum_{ \begin{array}{c} \scriptstyle 1 \leq b_1 \leq b_2 \leq \cdots \leq b_l \\
\scriptstyle a_i
> a_{i+1} \Rightarrow b_{i+1} > b_i \end{array}} x_{b_1} x_{b_2} \cdots
x_{b_l}.
\]
We have indexed the $F_{w^{-1}}(X)$ by the inverse permutation to
agree with the definition we shall give later.  Note that the
length $l(w)$ is equal to the degree of $F_w$ and the number
$|R(w)|$ of reduced decompositions of $w$ is given by the
coefficient of $x_1x_2\cdots x_l$ in $F_w$.

\begin{thm}[\cite{Sta84}]
\label{thm:Sta} The following properties of the generating
function $F_{w}$ hold for each $w \in S_n$:
\begin{enumerate}
\item
$F_w(X)$ is a symmetric function in $(x_1,x_2, \ldots)$.
\item
\label{it:2} Define $a_{w\ll} \in \Z$ by $F_w(X) = \sum_\ll
a_{w\ll}s_\ll(X)$. Then there exists partitions $\ll(w)$ and
$\mu(w)$ so that $a_{w\ll(w)} = a_{w\mu(w)} = 1$ and
\[
F_w(X) = \sum_{\ll(w) \preceq \ll \preceq \mu(w)} a_{w\ll}
s_\ll(X).
\]
\item
\label{it:3} Define an involution $*:S_n \to S_n$ by
$*:w_1w_2\cdots w_n \mapsto
(n+1-w_n)(n+1-w_{n-1})\cdots(n+1-w_1)$. Then
\[
\omega(F_w) = F_{w^*}.
\]
\item \label{it:4} We have
\[
s_1^\perp \cdot F_w = \sum_{w \gtrdot v} F_v.
\]
\end{enumerate}
\end{thm}

Edelman and Greene and separately Lascoux and Sch\"{u}tzenberger
showed the following (significantly harder) result concerning the
coefficients $a_{w\ll}$.
\begin{thm}[\cite{EG} and \cite{LS}]\label{thm:pos}
The coefficients $a_{w\ll}$ are non-negative.
\end{thm}

We now give a different formulation of the definition in a manner
similar to \cite{FS}.  Let $\C[S_n]$ denote the group algebra of
the symmetric group equipped with a inner product $\ip{w,v} =
\delta_{wv}$. Define linear operators $u_i: \C[S_n] \rightarrow
\C[S_n]$ for $i \in [1,n-1]$ by
\[
u_i.w = \begin{cases} s_i.w & \mbox{if $l(s_i.w) > l(w)$,} \\ 0 &
\mbox{otherwise.} \end{cases}
\]
The operators satisfy the braid relations $u_i u_{i+1} u_i =
u_{i+1} u_i u_{i+1}$ together with $u_i^2 = 0$ and $u_i u_j = u_j
u_i$ for $|i-j| \geq 2$.  They generate an algebra known as the
\emph{nilCoxeter algebra}.  Note that the action on $\C[S_n]$ is a
faithful representation of these relations.

Let $A_k(u) = \sum_{b_1 > b_2 > \cdots > b_k} u_{b_1} u_{b_2}
\cdots u_{b_k}$.  Then the Stanley symmetric functions can be
written as
\begin{equation}
\label{eq:def2} F_{w}(X) = \sum_{a = (a_1,a_2,\ldots,a_t)}
\ip{A_{a_t}(u)A_{a_{t-1}}(u) \cdots A_{a_1}(u) \cdot 1, w}
x_1^{a_1} x_2^{a_2} \cdots x_t^{a_t}
\end{equation}
where the sum is over all compositions $a$.  The symmetry of
$F_w(X)$ is then a consequence of the fact that the $A_k(u)$
commute.

For completeness, we explain briefly the relationship between
$F_w(X)$ and the Schubert polynomials of Lascoux and
Sch\"{u}tzenberger (see \cite{BJS}).  For $w \in S_n$, we have a
Schubert polynomial $\mathfrak{S}_w \in
\C[x_1,x_2,\ldots,x_{n-1}]$.  If $w \in S_n$, then $w \times 1^s
\in S_{n+s}$ denotes the corresponding permutation of $S_{n+s}$
acting on the elements $[1,n]$ of $[1,n+s]$. Similarly, $1^s
\times w \in S_{n+s}$ denotes the corresponding permutation acting
on the elements $[s+1,n+s]$ of $[1,n+s]$.  Schubert polynomials
have the important stability property $\mathfrak{S}_w =
\mathfrak{S}_{w \times 1^s}$. Stanley symmetric functions $F_w(X)$
are obtained by taking the other limit: $F_{w} = \lim_{s
\rightarrow \infty} \mathfrak{S}_{1^s \times w}$.  The limit is
taken by treating both sides as formal power series and taking the
limit of each coefficient.

\section{Affine Stanley symmetric functions}
\label{sec:affineStan}
Our first definition of affine Stanley symmetric functions will
imitate the definition (\ref{eq:def2}) above.  Let $\U_n$ be the
\emph{affine nilCoxeter algebra} generated over $\C$ by generators
$u_0,u_1,\ldots,u_{n-1}$ satisfying
\begin{align*}
u_i^2 & = 0 & \mbox{for all $i \in [0,n-1]$,} \\
u_i u_{i+1} u_i & = u_{i+1}u_iu_{i+1} &  \mbox{for all $i \in [0,n-1]$,} \\
u_i u_j &= u_j u_i & \mbox{for all $i,j \in [0,n-1]$ satisfying
$|i-j| \geq 2$.} \end{align*} Here and henceforth the indices are
to be taken modulo $n$.  A basis of $\U_n$ is given by the
elements $u_w = u_{\rho_1}u_{\rho_2}\cdots u_{\rho_l}$ where $\rho
= (\rho_1 \rho_2 \cdots \rho_l)$ is some reduced word for $w$ (see
\cite[Chapter 7]{Hum}).  The element $u_w \in \U_n$ does not
depend on the choice of reduced word $\rho$.

Let $a = a_1a_2\cdots a_k$ be a word with letters from $[0,n-1]$
so that $a_i \neq a_j$ for $i \neq j$.  Let $A =
\set{a_1,a_2,\ldots, a_k} \subset [0,n-1]$.  The word $a$ is
\emph{cyclically decreasing} if for every $i$ such that $i,i+1 \in
A$, the letter $i+1$ precedes $i$ in $a$.  We will call an element
$u \in \U_n$ cyclically decreasing if $u = u_a = u_{a_1}\cdots
u_{a_k}$ for some cyclically decreasing word $a$.  If $u$ is
cyclically decreasing and $u = u_{a_1}\cdots u_{a_k}$ then
necessarily $a = a_1 \cdots a_k$ will be cyclically decreasing.
The element $u$ is completely determined by the set $A =
\set{a_1,a_2,\ldots, a_k} \subset [0,n-1]$ and we write $u = u_A$.
Replacing $u_i$ by $s_i$ we make similar definitions of cyclically
decreasing affine permutations for the affine symmetric group.

Define $h_k(\u) \in \U_n$ for $k \in [0,n-1]$ by
\[
h_k(\u) = \sum_{A \in {[0,n-1] \choose k}} u_A
\]
where the sum is over subsets of $[0,n-1]$ of size $k$.  For
example if $n = 9$ and $A = \set{0,2,4,5,6,8}$ then $u_A = u_0 u_8
u_2 u_6 u_5 u_4 = u_2 u_6 u_5 u_4 u_0 u_8 = \cdots$.  A related
formula was given in \cite{Pos}, in the context of the
\emph{affine nil-Temperley-Lieb algebra}.  The affine
nil-Temperley-Lieb algebra is a quotient of the affine nilCoxeter
algebra given by the additional relations $u_iu_{i+1}u_i =
u_{i+1}u_iu_{i+1} = 0$.

Define a representation of $\U_n$ on $\C[\tS_n]$ by
\[
u_i.w = \begin{cases} s_i. w& \mbox{if $l(s_i.w) > l(w)$,} \\
0 & \mbox{otherwise.} \end{cases}
\]
It is easy to see that this is indeed a representation of $\U_n$.
If we identify $w \in \C[\tS_n]$ with $u_w \in \U_n$ then this is
essentially the left regular representation of $\U_n$.  Equip
$\C[\tS_n]$ with the inner product $\ip{w,v} = \delta_{wv}$. The
following definition was heavily influenced by \cite{FG}.

\begin{definition}
\label{def:affineStanley} Let $w \in \tS_n$.  Define the
\emph{affine Stanley symmetric functions} $\tF_w(X)$ by
\[\tF_{w}(X) = \sum_{a = (a_1,a_2,\ldots,a_t)}
\ip{h_{a_t}(u)h_{a_{t-1}}(u) \cdots h_{a_1}(u) \cdot 1, w}
x_1^{a_1} x_2^{a_2} \cdots x_t^{a_t},\] where the sum is over
compositions of $l(w)$ satisfying $a_i \in [0,n-1]$.
\end{definition}
The seemingly more general ``skew'' affine Stanley symmetric
functions \[\tF_{w/v}(X) = \sum_{a = (a_1,a_2,\ldots,a_t)}
\ip{h_{a_t}(u)h_{a_{t-1}}(u) \cdots h_{a_1}(u) \cdot v, w}
x_1^{a_1} x_2^{a_2} \cdots x_t^{a_t}\] are actually equal to the
usual affine Stanley symmetric functions $\tF_{wv^{-1}}(X)$.

Two properties follow straight from the definition.
\begin{prop}
Let $w \in \tS_n$.  Then the coefficient of $x_1x_2\cdots
x_{l(w)}$ in $\tF_w(X)$ is equal to the number of reduced words of
$w$.
\end{prop}

\begin{prop}
Suppose $w \in \Sn \subset \tS_n$.  Then $\tF_{w}(X) = F_{w}(X)$.
\end{prop}

The main theorem of this section is the following.
\begin{thm}
\label{thm:sym} The generating functions $\tF_w(X) \in \Ln$ are
symmetric.
\end{thm}

Theorem \ref{thm:sym} follows immediately from Proposition
\ref{prop:commute}.  In the following, intervals $[a,b]$ are to be
taken in the cyclic fashion within $[0,n-1]$.  Also, $\max$ and
$\min$ of a cyclic interval is meant to be taken modulo $n$ in the
obvious manner. So if $n = 6$ then $[4,1] = \set{4,5,0,1}$ and
$\max([4,1]) = 1$ and $\min([4,1]) = 4$.  We will need a technical
lemma first.

\begin{lem}
\label{lem:reduced} We have the following identities for reduced
words. \begin{enumerate} \item Let $a, b \in [0,n-1]$ with $a \neq
b-1$.  Then \[ a(a-1)(a-2)\cdots b a (a-1)(a-2) \cdots b \sim 0.\]
\item Let $a, b, c \in [0,n-1]$ satisfying $a \neq b-1$; $c \neq b$ and $c \in
[b,a]$.  Then \[a (a-1) (a-2) \cdots b c \sim (c-1)a (a-1) (a-2)
\cdots b.\]
\end{enumerate}
\end{lem}
\begin{proof}
Both results can be calculated by induction.
\end{proof}
So for example, the element $(s_4s_3s_2)(s_4s_3s_2)$ is not
reduced and we have $(s_6s_5s_4s_3s_2s_1)s_4 = s_3
(s_6s_5s_4s_3s_2s_1)$.

\begin{prop}
\label{prop:commute} The elements $h_k(\u)$ for $k \in [0,n -1]$
commute.
\end{prop}
\begin{proof}
For each $w \in \tS_n$ satisfying $l(w) = x+y$, we calculate the
coefficient of $u_w$ in $h_x(\u)h_y(\u)$ and $h_y(\u)h_x(\u)$.  We
assume that $x$ and $y$ are both not equal to $0$ for otherwise
the result is obvious. Let $u_w = u_A u_B$ where $|A| = x$ and
$|B| = y$. We need to exhibit a bijection between reduced
decompositions of this form and those of the form $u_w = u_Cu_D$
with $|C| = y$ and $|D| = x$.
We assume for simplicity (though it is not crucial to our proof)
that $A \cup B = [0,n-1]$ for otherwise we are in the non-affine
case and the proposition follows from results of Stanley
\cite{Sta84} or Fomin-Greene \cite{FG}. Let $A = \bigcup_i A_i$
and $B = \bigcup_i B_i$ be minimal decompositions of $A$ and $B$
into cyclic intervals.  If $A_i \subset B_j$ for some pair $(i,j)$
then we call $A_i$ an \emph{inner} interval and similarly for $B_k
\subset A_l$. Otherwise the interval is called \emph{outer}.

Using Lemma \ref{lem:reduced} and our assumption that $A \cup B =
[0,n-1]$ we can describe the outer intervals in an explicit
manner.  Each outer interval $A_i$ \emph{touches} an outer
interval $\rn(A_i)  = B_k$ called the \emph{right neighbour} of
$A_i$, for a unique $k$, so that $\min(A_i) = \max(B_k) +1$.  Also
$A_i$ overlaps with an outer interval $\ln(A_i) = B_l $ for a
unique $l$, so that $\max(A_i) \geq \min(B_l) - 1$ called the
\emph{left neighbour}.  If $\rn(A_i) = B_k$ then we also write
$A_i = \ln(B_k)$ and similarly for $\rn(B_k)$. Note that it is
possible that $\rn(A_i) = \ln(A_i)$ since we are working
cyclically.

Our bijection will depend only locally on each pair of an outer
interval $A^*$ and its right neighbour $B^* = \rn(A^*)$.  We call
the interval $I = [\min(B^*), \min(\ln(A^*)) - 1]$ a
\emph{critical interval}.  Critical intervals cover $[0,n-1]$ in a
disjoint manner.  For example, suppose $n = 10$ and $A =
\set{1,2,3,6,7,8,9}$ and $B = \set{0,2,4,5,7,9}$ (Figure
\ref{fig:intervals}), so that $u_A u_B = u_9u_8u_7u_6u_3u_2u_1
u_0u_9 u_7 u_5 u_4 u_2 u_0$. Then $A_1 = [1,3]$ and $A_2 = [6,9]$
are both outer intervals. Also $B_1 = [2,5]$, $B_2 = \set{7}$ and
$B_3 = [9,0]$. Only $B_2$ is an inner interval. The left neighbour
of $A_1$ is $\ln(A_1) = B_1$ and the right neighbour is $\rn(A_1)
= B_3$.  The critical intervals are $[9,1]$ and $[2,8]$.

Let $a = \min(\ln(A^*)) - 1$ and $b = \min(B^*)$. Let $c =
|[b,a]|$, $d = |A \cap [b,a]|$ and $e = |B \cap[b,a]|$. Renaming
for convenience, we let $S_1,S_2,\ldots,S_r$ be the inner
intervals (of $B$) contained in $A^*$ and $T_1,\ldots,T_t$ be
those contained in $B^*$, arranged so that $S_k > S_{k+1}$ for all
$k$ within $[b,a]$ and similarly $T_k > T_{k+1}$. We now define a
subset $U \subset [b,a]$ satisfying $|U| = d$. The algorithm
begins with $U = [b,a]$ and a changing index $i$ set to $i := a$
to begin with. The index $i$ decreases from $a$ to $b$ and at each
step the element $i$ may be removed from $U$ according to the
rule: \begin{enumerate} \item If $i \in A^*$ then we remove it
from $U$ unless $i \in S_k$ for some $k \in [1,r]$.
\item If $i \in B^*$ then we remove it from $U$ unless $i \in
T_k + 1$ for some $k \in [1,t]$.
\item Otherwise we do not remove $i$ from $U$ and set $i:= i-1$.
Repeat.
\end{enumerate}
When $|U| = d$ we stop the algorithm.  The algorithm always
terminates with $|U| = d$ since there are at least $c-d = |[b,a]|
- (A \cap [b,a])$ elements to remove.  In fact the algorithm
terminates before $i = b$ since $\cup_i S_i \neq A^* \cap I$.  We
will denote the result of the algorithm by $\phi(A^* \cup_i
T_i,B^* \cup_i S_i) := U$.  Note that $\min(U) = b$.

The bijection $u_Au_B \mapsto u_Cu_D$ is obtained by letting $D
\subset [0,n - 1]$ be the subset obtained from $B$ by changing $B
\cap I$ in each critical interval $I$ to $U$.  By the definition
of $U$ we see that $|D| = |A|$.  We claim that $u_Au_B = u_C u_D$
or alternatively $s_As_B(s_D)^{-1} = s_C$ for some $C$ satisfying
$|C| = |B|$ (here it is slightly more convenient to calculate
within the affine symmetric group, which is legal since our words
are all reduced). We can calculate this locally on each critical
interval since the $s_{D \cap I}$ commute as $I$ varies over
critical intervals. Note that $U$ always has the form of a
disjoint union $S_1 \cup S_2 \cup \cdots \cup S_{r'} \cup [b,a']$
for some $r' \leq r$ where $a'
> \max(B^*)$ or the form $S_1 \cup \cdots \cup S_r \cup \set{T_1 +
1} \cup \set{T_2+1} \cup \cdots \cup \set{T_{t'} +1 } \cup [b,a']$
where $a' \leq \max(B^*)$.

Let us assume that $U$ has the first form.  Focusing on $I = [b,a]
=  [\min(B^*), \min(\ln(A^*)) - 1] $ we are interested in
\[\underline{s} = s_{A^* \cap I} s_{T_1} \cdots s_{T_t} s_{S_1}
\cdots s_{S_r} s_{B^*} (s_{[b,a']})^{-1} (s_{S_{r'}})^{-1} \cdots
(s_{S_1})^{-1}.\] Then we get
\begin{align*} \underline{s} &= s_{A^* \cap I} s_{T_1} \cdots s_{T_t} s_{S_{r' + 1}}
\cdots s_{S_r} (s_{[\max(B^*)+1,a']})^{-1} \\ &=  s_{S_{r' + 1} -
1} \cdots s_{S_r - 1 } s_{T_1} \cdots s_{T_t} s_{A^* \cap I}
(s_{[\max(B^*)+1,a']})^{-1}
\\ &= s_{S_{r'+ 1} - 1} \cdots s_{S_r - 1 } s_{T_1} \cdots s_{T_t}
s_{[a' + 1, a ]} & \mbox{using $\max(B^*) +1 = \min(A^*)$}.
\end{align*}
We used Lemma \ref{lem:reduced} repeatedly and also the fact that
the certain intervals do not ``touch'' and so commute. Let $U'$ be
the disjoint union $[a' + 1, a ] \cup \set{S_{s' + 1} - 1} \cup
\cdots \cup \set{S_s - 1} \cup T_1 \cup \cdots T_t$.  Note that it
is always the case that $\max(U') = a$.  The other form of $U$
involves a similar calculation. One checks that we can combine
this argument for each critical interval showing that
$s_As_B(s_D)^{-1}$ is indeed equal to $s_C$ for some $C$.

Finally, we need to show that this map is a bijection.  Again we
work locally on a critical interval and assume that $U$ has the
first form. If we replace $A^*$ (more precisely $A^* \cap I$) by
$U'$ and $B^*$ by $U$, then our internal intervals are $S'_1 =
S_1, \ldots S'_{r'} = S_{r'}$ and $T'_1 = S_{r'+1} - 1, \ldots,
T'_{r-r'} = S_{r'}  -1, T'_{r-r'+1} = T_1, \ldots T'_{r-r'+t} =
T_t$.  We now show that $B^* \cup_i S_i = \phi(U',U)$ from which
the bijectivity will follow.  Note that since $\min(U) = b$ and
$\max(U') = a$ the critical intervals of $u_Cu_D$ are the same as
those of $u_Au_B$.  By definition $\phi(U',U)$ keeps
$S'_1,S'_2,\ldots$ and keeps $T'_1 +1,T'_2-1,\ldots, T'_{r-r'} +
1$, removing all other values up to this point.  At this point the
algorithm stops since $\phi(U',U)$ is of the correct size.  We see
that we obtain $\phi(U',U) = B^* \cup_i S_i$ back in this way. A
similar argument works for the second form of $U$.
\end{proof}

\begin{figure}[ht]
\pspicture(-50,-50)(50,50)
 \pscircle[linewidth=0.5pt](0,0){40}  \SpecialCoor
\rput(50;0){0} \rput(50;324){1} \rput(50;288){2} \rput(50;252){3}
\rput(50;216){4} \rput(50;180){5} \rput(50;144){6}
\rput(50;108){7} \rput(50;72){8} \rput(50;36){9}

\psdot[dotscale=1.5,dotstyle=square](40;36)
\psdot[dotscale=1.5,dotstyle=square](40;0)
\psdot[dotscale=1.5,dotstyle=square](40;288)
\psdot[dotscale=1.5,dotstyle=square](40;252)
\psdot[dotscale=1.5,dotstyle=square](40;216)
\psdot[dotscale=1.5,dotstyle=square](40;180)
\psdot[dotscale=1.5,dotstyle=square](40;108)

\psdot[dotscale=1](40;36) \psdot[dotscale=1](40;72)
\psdot[dotscale=1](40;108)\psdot[dotscale=1](40;144)
\psdot[dotscale=1](40;252)\psdot[dotscale=1](40;288)
\psdot[dotscale=1](40;324)

\endpspicture
\caption{Dots represent elements of $A$.  Squares represent
elements of $B$.} \label{fig:intervals}
\end{figure}

\begin{example}
We illustrate the map $U = \phi(A^* \cup_i T_i,B^* \cup_i S_i)$ of
the proof. Suppose $[b,a] = [2,20]$ and $A^* = [14,20]$, $B^* =
[2,13]$.  Let $S_1 = [16,18]$ and $T_1 = [8,11]$ and $T_2 =
\set{5}$ be the inner intervals.  Then $d = 12$ and $U =
\set{2,3,4,5,6,9,10,11,12,16,17,18}$.  We can compute that
\[
s_{A^*}s_{11}s_{10}s_9 s_8 s_5 s_{B^*} s_{18}s_{17}s_{16} s_2 s_3
s_4 s_5 s_6 s_9 s_{10} s_{11} s_{12} s_{16} s_{17} s_{18} =
s_{A^*}s_{[7,13]}s_5
\]
so that $U' = [7,20] \cup \set{5}$.  Finally one checks that $B^*
\cup_i S_i = \phi(U',U)$.
\end{example}

We end this section by giving two alternative descriptions of the
affine Stanley symmetric functions, the first one imitating the
original definition of Stanley. Let $w \in \tS_n$. Let $a =
(a_1,\ldots,a_l) \in R(w)$ be a reduced word and $b = (b_1 \geq
b_2 \cdots \geq b_l)$ be an positive integer sequence. Then
$(a,b)$ is called a \emph{compatible pair for $w$} if whenever
$b_i = b_{i+1} = \cdots = b_j$ and $\set{k,k+1} \subset
\set{a_i,a_{i+1},\ldots,a_j}$ then we have that $k+1$ precedes $k$
(for any $i,j,k$). Two compatible pairs $(a,b)$ and $(a',b')$ are
equivalent if $b = b'$ and for any maximal interval $[i,j] \subset
[1,l]$ satisfying $b_i = b_{i+1} = \cdots = b_j$ we have that $a_i
a_{i+1} \cdots a_j$ and $a'_i a'_{i+1} \cdots a'_j$ are reduced
words for the same affine permutation.  The following proposition
is clear from the definitions.
\begin{prop}[Alternative Definition 1]  The affine Stanley
symmetric functions are given by
\[
\tF_w(X) = \sum_{ \overline{ (a,b)}} x_{b_1} x_{b_2} \cdots
x_{b_l}
\]
where the sum is over equivalence classes $\overline{ (a,b)}$ of
compatible pairs for $w$.
\end{prop}

Now let $w\in \tS_n$ of length $l$ and suppose $\a =
(\a_1,\a_2,\ldots,\a_r)$ is a composition of $l$.  An
\emph{$\a$-decomposition} of $w$ is an ordered $r$-tuple of
cyclically decreasing affine permutations $(w^1,w^2,\ldots,w^r)\in
\tS_n^r$ satisfying $l(w^i) = \a_i$ and $w = w^1w^2\cdots w^r$.
The following alternative definition is also immediate.
\begin{prop}[Alternative Definition 2] The affine Stanley
symmetric function $\tF_w(X)$ is given by
\[
\tF_w(X) = \sum_{\a} (\mbox{number of $\a$-decompositions of
$w$})\cdot x^\a
\]
where the sum is over all compositions $\a$ of $l$.
\end{prop}

\section{Representations of the affine nilCoxeter algebra}
\label{sec:Vrep} Let $V$ be a complex representation of $\U_n$
with a distinguished basis $\{v_p \mid p \in P\}$ for some
indexing set $P$. Let $\ip{.,.}: V \times V \rightarrow \C$ be the
inner product defined by $\ip{v_p,v_q} = \delta_{pq}$ for $p,q \in
P$. For any $p, q \in P$ one can define \emph{$V$-affine Stanley
symmetric functions} by $\tF_{q/p}(X) \in \Ln$ by
\[\tF_{q/p}(X) = \sum_{a = (a_1,a_2,\ldots,a_t)}
\ip{h_{a_t}(u)h_{a_{t-1}}(u) \cdots h_{a_1}(u) \cdot v_p, v_q}
x_1^{a_1} x_2^{a_2} \cdots x_t^{a_t},\] where the sum is over
compositions of $l(w)$ satisfying $a_i \in [0,n-1]$.  By
Proposition \ref{prop:commute} these functions are indeed
symmetric functions.

\begin{prop}
\label{prop:Vaffine} Suppose $u_w \cdot v_p = v_q$ and $w \in
\tS_n$ is the only affine permutation such that $\ip{u_w\cdot
v_p,v_q} \neq 0$. Then $\tF_{q/p}(X) = \tF_w(X)$.
\end{prop}
\begin{proof}
For each composition $a = (a_1,a_2,\ldots,a_t)$ expand
$h_{a_t}(u)h_{a_{t-1}}(u) \cdots h_{a_1}(u)$ in the basis
$\{u_v\}$ of $\U_n$.  Using the assumption, the proposition
follows immediately upon comparison with Definition
\ref{def:affineStanley}.
\end{proof}

More generally, for arbitrary $v_p$ and $v_q$ let $c_w = \ip{u_w
\cdot v_p, v_q}$. Then $\tF_{q/p} = \sum_{w \in \tS_n} c_w\tF_w$.
We have not found any interesting generating functions of this
form.

If in addition $u_i$ acts on the basis $\{v_p\}_{p\in P}$ with
non-negative matrix coefficients, then $\tF_{p/q}$ will be
monomial-positive.  This will be the case for all the
representations of $\U_n$ that we will be considering.

\section{Coproduct}
\label{sec:basic} We now give the analogue of part (\ref{it:4}) of
Theorem \ref{thm:Sta}.

\begin{thm}[Coproduct formula]
The following coproduct expansion holds:
\[ \tF_w(x_1,x_2,\ldots,y_1,y_2,\ldots) = \sum_{uv = w} \tF_v(x_1,x_2,\ldots)\tF_u(y_1,y_2,\ldots).\]
In particular we have
\[
s_1^\perp \tF_w = \sum_{w \gtrdot v} \tF_v.
\]
\end{thm}
\begin{proof}
The first formula follows immediately from the definition and the
fact that $\tF_{w/v}(Y) = \tF_{wv^{-1}}(Y)$.  To obtain the second
formula, we first write, using the first formula and
(\ref{eq:coproduct}),
\[
\sum_{uv = w} \tF_v(X)\otimes\tF_u(Y) = \sum_\ll s_\ll^\perp(X)
\tF_w(X) \otimes s_\ll(Y).
\]
The terms of the formula are to be interpreted within $\Lambda$,
even though the sum is an element of $\Ln$.  Now take the inner
product of both sides with $s_1(Y)$ to get
\[
 s_1^\perp(X) \tF_w(X) = \sum_{uv = w} \tF_v(X)\ip{\tF_u(Y),s_1(Y)}.
\]
Now $\ip{\tF_u(Y),s_1(Y)} = 0$ unless $u = s_i$ is a simple
reflection for some $i$, in which case $\tF_{s_i}(Y) = s_1(Y)$.
This gives the second formula.
\end{proof}

\section{Monomial dominance}
\label{sec:monomialDominance} We now show that there is a dominant
term in the \emph{monomial} expansion of an affine Stanley
symmetric function $\tF_w(X)$.  Let $c'(w)=c(w^{-1})$ denote the
code of the inverse $w^{-1}$ of $w$, so that $c'_{w(i)} =
\#\set{j:\; j < i \; \mbox{and} \; w(j) > w(i)}$. Let $\mu(w)$
denote the partition which is conjugate to the decreasing
permutation of $c'(w)$.

\begin{thm}
\label{thm:monomial} Let $w \in \tS_n$.  Then
\begin{enumerate}
\item
If $[m_\ll]\tF_w \neq 0$ then $\ll \preceq \mu(w)$.
\item
We have $ [m_\mu(w)]\tF_w = 1$.
\end{enumerate}
\end{thm}

\begin{proof}
Left multiplication of $w$ by $s_i$ acts on $c'(w)$ by
\[
s_i: (c'_1,\ldots,c'_{i},c'_{i+1},\ldots,c'_{n})\longmapsto
(c'_1,\ldots,c'_{i+1}+1,c'_{i},\ldots,c'_{n})
\]
whenever $l(s_iw) > l(w)$.  Applying a term of $h_k(\u)$ to $w$
will increase $k$ different entries of $c'(w)$ by 1 and also
permute the entries (assuming the result is non-zero), since $u_i$
never acts after $u_{i+1}$. Using this repeatedly we see that if
$m_\ll$ occurs in $\tF_w$, we must have $\mu_1(w) \geq \ll_1$ and
then $\mu_1(w) + \mu_2(w) \geq \ll_1 + \ll_2$ and so on.  So $\ll
\preceq \mu(w)$.

Now we check that the coefficient of $x^{\mu(w)}$ in $\tF_w(X)$ is
1.  To see this, we work by going down in the Bruhat order or
equivalently, acting on $w$ by $h_k(\u)^\perp$ (the adjoint with
respect to $\ip{.,.}$ of $h_k(\u)$). Multiplying $w$ by a term of
$h_{\mu_1(w)}(\u)^\perp$ means decreasing $\mu_1(w)$ different
entries of $c'(w)$ by 1 each (and also permuting the entries in
some way). But $c'(w)$ only has $\mu_1(w)$ non-zero entries, and
so there is only one possible resulting code $c'(v)$: it is
obtained from $c'(w)$ by taking all non-zero entries $c'_i$ and
shifting them each to the right (cyclically) one entry.  This is
because entries can only decrease (by 1) by shifting to the right,
and once such an entry is shifted we are forbidding it from moving
again. Now the conjugate of the decreasing permutation of $c'(v)$
is exactly $(\mu_2(w),\mu_3(w),\ldots)$ so our result follows from
induction.
\end{proof}

\begin{corollary}
\label{cor:nL} The subalgebra $\nL(\u)$ generated of $\U_n$ by
$\{h_k(\u)\}_{k = 1}^{n-1}$ is isomorphic to $\nL$ with
isomorphism given by $h_i \mapsto h_i(\u)$ for $1 \leq i \leq
n-1$.
\end{corollary}
\begin{proof}
Suppose to the contrary that the $h_k(\u)$ are not algebraically
independent. Then there is some relation $h_\ll(\u) = \sum_\nu
a_\nu h_\nu(\u)$ where we may pick $\ll$ so that no $\nu$
appearing on the right hand side satisfies $\nu \prec \ll$.  Now
pick $w$ so that $\mu(w) = \ll$.  Then by Theorem
\ref{thm:monomial}, $1 = \ip{h_\ll(\u) \cdot 1, w} = \sum_\nu
a_\nu \ip{h_\nu(\u) \cdot 1, w} = 0$, a contradiction.
\end{proof}

We denote by $f(\u)$ the image of $f \in \nL$ under the
isomorphism $\nL \cong \nL(\u)$.

\section{Conjugacy}
\label{sec:conjugacy}

Define $\omega: \nL \rightarrow \nL$ as usual by $\omega: h_i
\mapsto e_i$.  Define $\omega^+: \Ln \rightarrow \Ln$ by requiring
that $\ip{\omega(f), \omega^+(g)} = \ip{f,g}$ where $f \in \nL$
and $g \in \Ln$.  Alternatively, we require that the sets
$\set{e_\ll \mid \ll \in \pa^n}$ and $\set{\omega^+(m_\ll) \mid
\ll \in\pa^n}$ form dual bases of $\nL$ and $\Ln$.  The map
$\omega^+$ is clearly an involution but it does not agree with
$\omega$ (see for example \cite[Chapter 7, Ex. 9]{EC2}).

Denote by $w \mapsto w^*$ the involution of $\tS_n$ given by $s_i
\mapsto s_{n-i}$ (with $s_0 \mapsto s_0$).  In terms of the window
realization of $\tS_n$, we have $[w(1),w(2),\ldots,w(n)]^* =
[n+1-w(n),n+1-w(n-1),\ldots,n+1-w(1)]$.  Similarly, $u_i \mapsto
u_{n-i}$ defines an algebra involution (also denoted $*$) of
$\U_n$.

\begin{thm}[Conjugacy formula]
\label{thm:conj} Let $w \in \tS_n$.  Then $\omega^+(\tF_w) =
\tF_{w^*}$.
\end{thm}

We shall prove Theorem \ref{thm:conj} by calculating within the
subalgebra $\nL(\u)$ of Corollary \ref{cor:nL}.  The following
result says that $e_k(\u) = (h_k(\u))^*$.

\begin{prop}
\label{prop:eformula} The elements $e_k(\u) \in \U_n$ are given by
\[ e_k(\u) = \sum_{A \in {[0,n-1] \choose k}} \tilde{u}_A,
\]
where for a $k$-subset $A = \set{a_1,a_2,\ldots,a_k} \subset
[0,n-1]$ the element $\tilde{u}_A \in \U_n$ is defined as any
expression $u_{a_1} u_{a_2} \cdots u_{a_k}$ where if $i$ and $i+1$
(modulo $n$) are both in $A$ then $u_{i}$ must precede $u_{i+1}$
within $\tilde{u}_A$.
\end{prop}

\begin{proof}
We verify this using the relation \begin{equation} \label{eq:eh}
e_k(\u) = h_k(\u) - h_{k-1}(\u)e_1(\u) + \cdots \pm
h_1(\u)e_{k-1}(\u).\end{equation} First, we restrict our attention
to the monomials which only involve the set of generators
$\set{u_1,u_2,\ldots,u_{n-1}}$.  Then one may write
\[h_k(\u) = \sum_{n-1 \geq \{a_1 > a_2 > \cdots > a_k\} \geq 1}
u_{a_1}u_{a_2} \cdots u_{a_k}\] and we assume that
\[e_l(\u) = \sum_{n-1 \geq \{a_1 < a_2 < \cdots < a_l\} \geq 1}
u_{a_1}u_{a_2} \cdots u_{a_l}\] is known for $l < k$.  (The base
case $k = 1$ is clear.)  Now for $k > l \geq 1$,
$h_{k-l}(\u)e_l(\u)$ can be written as $A_l + B_l$ where
\[
A_l = \sum_{n-1 \geq \set{a_1 > a_2 > \cdots > a_{k-l} < a_{k-l+1}
< \cdots < a_k } \geq 1} u_{a_1}u_{a_2} \cdots u_{a_k}
\]
and
\[
B_l = \sum_{n-1 \geq \set{ a_1 > a_2 > \cdots > a_{k-l} >
a_{k-l+1} < \cdots < a_k} \geq 1} u_{a_1}u_{a_2} \cdots u_{a_k}.
\]
Note that $h_k(\u)  = B_1$ and for $k > l \geq 1$, we have $A_l =
B_{l+1}$ so all but one of the terms on the right hand side of
(\ref{eq:eh}) cancel to give $e_k(\u) = A_{k-1}$, which is the
desired formula.  This proves the theorem when the monomials are
restricted to $\set{u_1,u_2,\ldots,u_{n-1}}$.  But since $k \leq
n-1$, any monomial $u_w$ in (\ref{eq:eh}) only involves a proper
subset of the generators $\set{u_0,u_1,u_2,\ldots,u_{n-1}}$, so we
can calculate the coefficient of that monomial in $e_k(\u)$ by
setting $u_i = 0$ for some $i$.  The theorem follows.
\end{proof}

More generally, when $\ll = (a,1^b)$ is a hook shape satisfying
$s_\ll \in \nL$ then $s_\ll(\u)$ can be written as a sum over the
reading words of certain tableaux (see \cite{Lam}). We shall not
need this generality; however, see Proposition \ref{prop:kschur}.

\begin{proof}[Proof of Theorem \ref{thm:conj}] Write the
\emph{affine non-commutative Cauchy kernel}
\[\Omega^{(n)}(x,\u):= \sum_{\ll \in \pa^n} h_\ll(\u) m_\ll(X) = \sum_{\ll \in \pa^n} e_\ll(\u)
\omega^+(m_\ll(X))\] where the second equality follows from the
definition of $e_\ll(\u)$ and an argument similar to \cite[Lemma
7.9.2]{EC2}.

By definition $\tF_w(X) = \ip{\Omega^{(n)}(x,\u) \cdot 1, w} =
\sum_{\ll \in \pa^n} \ip{e_\ll(\u) \cdot 1, w}\omega^+(m_\ll(X))$.
By Theorem \ref{prop:eformula}, $e_\ll(\u)$ is obtained from
$h_\ll(\u)$ by the involution $u_i \mapsto u_{n-1-i}$, so
$\ip{e_\ll(\u) \cdot 1, w} = \ip{h_\ll(\u) \cdot 1, w^*}$.  This
completes the proof of the theorem.
\end{proof}

For later use, we have the following proposition.
\begin{prop}
\label{prop:inverseStar} Let $w \in \tS_n$.  Then $\tF_{w^*} =
\tF_{w^{-1}}$.
\end{prop}
\begin{proof}
The reduced words of $w^{-1}$ are obtained by reversing the
reduced words of $w$.  But each term of $e_\ll(\u)$ is also
obtained from a term of $h_\ll(\u)$ by reversing the order of the
generators.  This shows that $\ip{h_\ll(\u)\cdot 1, w^{-1}} =
\ip{e_\ll(\u) \cdot 1, w} = \ip{h_\ll(\u) \cdot 1, w^*}.$
\end{proof}

Let $\Z/n\Z$ act on $\tS_n$ by the action $p.s_i = s_{i+p}$ for $p
\in \Z/n\Z$.  Since the definition of $h_k(\u)$ is invariant under
the analogous transformations of $\U_n$, we have the following
symmetry of affine Stanley symmetric functions:
\begin{prop}
Let $w \in \tS_n$ and $p \in \Z/n\Z$.  Then $\tF_w = \tF_{p.w}$.
\end{prop}

\section{Affine Schur functions}
\label{sec:affineSchur} A permutation $w \in \tS_n$ is
\emph{Grassmannian} (or more precisely \emph{left-Grassmannian})
if it is a minimal length coset representative for a coset of
$\Sn\backslash\tS_n$ where $\Sn \cong S_n$ is the maximal
parabolic subgroup generated by the $n-1$ generators
$s_1,\ldots,s_{n-1}$. By general facts concerning parabolic
subgroups of Coxeter groups \cite{Hum}, the minimal length coset
representative $\bar w$ of a coset $(\Sn) w$ is unique and
satisfies $l(u \bar w) = l(u) + l(\bar w)$ for any $u \in \Sn$.
There is a natural correspondence between the minimal length coset
representatives corresponding to another embedding of $S_n$ into
$\tS_n$, and the ones we have called Grassmannian.  In particular
the associated affine Stanley symmetric functions are equal under
this correspondence so we will only consider the Grassmannian
permutations.

A permutation $w$ is Grassmannian if left multiplication by $s_i$
always increases the length $l(w)$.  This is equivalent to $c'(w)$
being a weakly increasing sequence, or equivalently, that the
window $[w^{-1}(1),w^{-1}(2),\ldots,w^{-1}(n)]$ of $w^{-1}$ is
increasing. In fact, the correspondence $w \leftrightarrow \mu(w)$
is a bijection between Grassmannian permutations and $\pa^n$ (see
\cite{BB}).

\begin{definition}
An affine Stanley symmetric function $\tF_w(X)$ is called an
\emph{affine Schur function} if $w$ is a Grassmannian permutation.
If $\mu = \mu(w)$, we write $\tF_{\mu}(X) := \tF_w(X)$.
\end{definition}

Affine Schur functions had earlier been defined by Lapointe and
Morse in a different manner, and were called \emph{dual $k$-Schur
functions}. We will see the origin of this name later.

\begin{thm}
The affine Schur functions $\set{\tF_{\mu} :\; \mu \in \pa^n}$
form a basis of $\Ln$.
\end{thm}
\begin{proof}
By Theorem \ref{thm:monomial}, $\tF_\mu = \sum_{\ll \preceq \mu}
b_{\mu\ll} m_\ll$ for some coefficients $b_{\mu\ll} \in \Z$
satisfying $b_{\mu\mu} = 1$.  Since the transition matrix between
$\{\tF_\mu\}$ and $\set{m_\ll}$ is uni-triangular, the theorem
follows.
\end{proof}
Now define $a_{w\ll} \in \Z$ by
\[
\tF_{w}(X) = \sum_{\ll \in \pa^n} a_{w\ll}\tF_{\ll}(X).
\]
The fact that the coefficients $a_{w\ll}$ are integers follows
from the fact that the transition matrix between $\{\tF_\mu\}$ and
$\set{m_\ll}$ is uni-triangular with integer coefficients,
together with the fact that the monomial expansion of $\tF_w$ has
integer coefficients. Let $\tilde{f}^\ll = [x_1x_2\cdots
x_{l(u)}]\tF_\ll(X)$ be the number of reduced decompositions of
the Grassmannian permutation $u$ satisfying $\mu(u) = \ll$.  Thus
for any $w \in \tS_n$ we have
\[
\#R(w) = \sum_\ll a_{w\ll} \tilde{f}^\ll.
\]
In fact we conjecture that $a_{w\ll} \geq 0$; see Section
\ref{sec:positivity}.  In the non-affine case, the numbers
$\tilde{f}^\ll$ are dimensions of irreducible representations of
the symmetric group and are given by the well-known hook length
formula; see \cite{EC2}.  It is unknown whether a closed formula
for $\tilde{f}^\ll$ exists in the affine case, though
$\tilde{f}^\ll$ does count the number of certain tableaux, known
as $k$-tableaux; see Section \ref{sec:ktableaux}.

Let $w$ be a Grassmannian permutation.  Then since the involution
$*: \tS_n \rightarrow \tS_n$ sends $\Sn$ to $\Sn$, the permutation
$w^*$ is also a Grassmannian permutation.  We thus obtain an
involution $*: \pa^n \rightarrow \pa^n$ given by requiring that
$\mu(w)^* = \mu(w^*)$ for Grassmannian permutations $w$. Combining
this with Theorem \ref{thm:conj} we obtain
\begin{equation}
\label{eq:omegaStar} \omega^+(\tF_\ll) = \tF_{\ll^*}.
\end{equation}
Let $v$ be a minimal coset representative of a right coset in
$\tS_n/\Sn$ (a \emph{right-Grassmannian} permutation).  Since $v$
is the inverse of some Grassmannian permutation, by Proposition
\ref{prop:inverseStar}, the associated affine Stanley symmetric
function $\tF_v$ is equal to an affine Schur function so in fact
we have lost no generality considering the left-Grassmannian
permutations instead of the right-Grassmannian permutations.

The involution $*$ on $\pa^n$ has been studied in a different form
in \cite{LLM} where it is called \emph{$k$-conjugation}. Define
the partial order $\prec^*$ on $\pa^n$ by $\ll \prec^* \mu$ if and
only if $\mu^* \prec \ll^*$.  The partial order $\prec^*$ is not
the same as $\prec$.  For example $(2,2)$ and $(2,1,1)$ are both
fixed points of $*$ for $n = 3$ (the author thanks J.~Morse for
this example).

Let $\ll(w) = \mu(w^{-1})^*$ (note that $c(w^{-1})$ and $c(w^*)$
are rearrangements of each other so that $\mu(w^{-1}) =
\mu(w^*)$).

\begin{thm}[Dominant Terms]
\label{thm:dominance} Let $w \in \tS_n$.  Then
\begin{enumerate}
\item
If $a_{w\ll} \neq 0$ then $\ll(w) \preceq^* \ll \preceq \mu(w)$.
\item
We have $a_{w\mu(w)} = a_{w\ll(w)} = 1$.
\end{enumerate}
\end{thm}
\begin{proof}
The statements involving $\mu(w)$ follow from Theorem
\ref{thm:monomial} and the comments earlier.  Applying this to
$w^{-1}$, we have $\tF_{w^{-1}}(X) = \tF_{\mu(w^{-1})} + \sum_{\ll
\prec \mu(w^{-1})} a_{w^{-1}\ll} \tF_{\ll}(X)$. Applying
$\omega^+$ to both sides and using Theorem \ref{thm:conj},
Proposition \ref{prop:inverseStar} and (\ref{eq:omegaStar}) we get
\[
\tF_{w}(X) = \tF_{\ll(w)} + \sum_{\ll(w) \prec^* \ll^*}
a_{w^{-1}\ll} \tF_{\ll^*}(X)
\]
which implies the other statements of the Theorem.
\end{proof}

We end this section with a question: for which $w \in \tS_n$ is
$\mu(w)^* = \mu(w^*)$?  Is it the same as the class of
permutations $w \in \tS_n$ such that $\tF_w$ is equal to an affine
Schur function?  See also Problem \ref{prob:vex}.

\section{Affine symmetric group and $n$-cores}
\label{sec:ncores}  We now describe an action of the affine
symmetric group on partitions.  Further details for the material
of this section can be found in \cite{vL, Las}.

A \emph{$n$-ribbon} is a connected skew shape $\ll/\mu$ of size
$n$ which contains no $2 \times 2$ square. A partition $\ll$ is an
\emph{$n$-core} if no $n$-ribbon $\ll/\mu$ can be removed from it
to obtain another partition $\mu$.  Let $\p^n$ denote the set of
$n$-cores.

\begin{figure}
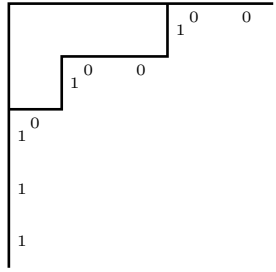

\pspicture(0,0)(100,100) \fontsize{6pt}{6pt} \rput(5,10){$1$}
\rput(5,30){$1$} \rput(5,50){$1$} \rput(10,55){$0$}
\rput(25,70){$1$} \rput(30,75){$0$} \rput(50,75){$0$}
\rput(65,90){$1$} \rput(70,95){$0$} \rput(90,95){$0$}
\psline(0,0)(0,100)(100,100)
\psline(0,60)(20,60)(20,80)(60,80)(60,100)
\endpspicture
\caption{The edge sequence $p(31) = (\ldots,
1,1,1,0,1,0,0,1,0,0,\ldots)$.} \label{fig:edgesequence}
\end{figure}

If $\ll$ is a partition, we let $p(\ll)$ denote the \emph{edge
sequence} of $\ll$.  The edge sequence $p(\ll) = (\ldots,
p_{-2},p_{-1},p_0,p_1,p_2,\ldots)$ is the doubly infinite bit
sequence obtained by drawing the partition in the English notation
and reading the ``edge'' of the partition from bottom left to top
right -- writing a 1 if you go up and writing a 0 if you go to the
right (see Figure \ref{fig:edgesequence}).  We shall normalise our
notation for edge sequences by requiring that the empty partition
$\emptyset$ has edge sequence $p(\emptyset)_i = 1$ for $i \leq 0$
and $p(\emptyset)_i = 0$ for $i \geq 1$.  Adding a box to a
partition corresponds to changing two adjacent entries of the edge
sequence $p_i,p_{i+1}$ from $(0,1)$ to $(1,0)$.
Adding a $n$-ribbon to a partition $\ll$ corresponds to finding an
index $i \in \Z$ such that $p_i(\ll) = 1$ and $p_{i+n}(\ll) = 0$,
then changing those two bits to $p_i(\ll) = 0$ and $p_{i+n}(\ll) =
1$.

Let $\ll$ be an $n$-core with edge sequence $p(\ll) = (\ldots,
p_{-2},p_{-1},p_0,p_1,p_2,\ldots)$.  Then there is no index $i$ so
that $p_i(\ll) = 0$ and $p_{i+n}(\ll) = 1$.  Equivalently, the
subsequences
\[p^{(i)}(\ll) = (\ldots,p_{i-2n},p_{i-n},p_i,p_{i+n},p_{i+2n},\ldots)\] all look
like $(\ldots,1,1,1,1,0,0,0,0,\ldots)$ with a suitable shift.
Define the \emph{offsets} $\{d_i = d_i(\ll) \mid i \in \Z\}$ by
requiring that $p_{i+nd_i} = 0$ and $p_{i+n(d_i-1)} = 1$.  The
offsets satisfy $d_{i-n} = d_i + 1$ and $d_1 + d_2 + \cdots +
d_{n} = 0$ and completely determine the $n$-core.

Now let $\p$ denote the set of doubly infinite $(0,1)$-sequences
$p = (\ldots, p_{-2},p_{-1},p_0,p_1,p_2,\ldots)$ and let $\C[\p]$
denote the space of formal $\C$- linear combinations of such
sequences.  Let $\tS_n$ act on $\p$ by letting $s_i$ act on $p =
(\ldots, p_{-2},p_{-1},p_0,p_1,p_2,\ldots)$ by swapping $p_{kn+i}$
and $p_{kn+i+1}$ for each $k \in \Z$.  One can check directly that
this defines a representation of $\tS_n$ on $\C[\p]$.

A sub-representation $\C[\p^*]$ of $\C[\p]$ is given by taking
only those bit sequences $p \in \p^*$ satisfying $p_N = 1$ for
sufficiently small $0 \gg N$ and $p_N = 0$ for $N \gg 0$.  These
sequences correspond to possibly shifted edge sequences of
partitions.  It is easy to see that $\C[\p^*]$ is indeed a
sub-representation, but it is by no means irreducible.  We let
$\tS_n$ act on partitions by the corresponding action on the edge
sequences. The action of $s_i \in \tS_n$ acts by adding and or
removing boxes along certain diagonals.

The proof of the following Proposition is straightforward.
\begin{prop}
\label{prop:ncores} The orbit $\tS_n \cdot \emptyset$ is equal to
the set of $n$-cores.  Let $\ll$ be an $n$-core with offsets
$d_i(\ll)$. Then $\mu = s_i \cdot \ll$ is an $n$-core with offsets
$d_j(\mu) = d_j(\ll)$ for $j \neq i,i+1$ and $d_{i+1}(\mu) =
d_i(\ll)$ and $d_i(\mu) = d_{i+1}(\ll)$.  If $d_i(\ll) >
d_{i+1}(\ll)$ then boxes are added; if $d_i(\ll) < d_{i+1}(\ll)$
then boxes are removed and if $d_i(\ll) = d_{i+1}(\ll)$ then $\ll
= \mu$.
\end{prop}

One can see (for example using Proposition \ref{prop:ncores}) that
the stabiliser of the empty partition is $\Sn \subset \tS_n$, so
the set $\p^n$ of $n$-cores is naturally isomorphic to
$\tS_n/\Sn$.  We may thus identify $n$-cores with
\emph{right-Grassmannian} permutations -- the set $\Snn$ of
minimal length coset representatives of $\tS_n/\Sn$.  If $w \in
\Snn$ satisfies $w \cdot \emptyset = \ll \in \p^n$ then we write
$w = w(\ll)$.

The following relation between the $n$-cores and the affine
symmetric group is known (see \cite{Las}).

\begin{prop}
\label{prop:Las} Let $\ll, \mu \in \p^n$ be $n$-cores.  Then $\ll
\subset \mu$ if and only if $w(\ll)$ is less than $w(\mu)$ in
(strong) Bruhat order.
\end{prop}

The action of $\tS_n$ on $\p^n$ corresponds to the left action of
$\tS_n$ on $\tS_n/\Sn$.  We will need the following general fact
for Coxeter groups.

\begin{lem}
\label{lem:cox} Let $W$ be a Coxeter group, $W_I$ a parabolic
subgroup and $W^I$ a the set of minimal length coset
representatives of $W/W_I$.  Let $w \in W^I$ and $s_i$ be a simple
generator.  Then either $s_i w \in W^I$ or $s_iw \in wW_I$.
\end{lem}
\begin{proof}
Let $l(w) = l$.  Suppose that $s_iw = vu$ for $v \in W^I$ and $u
\in W_I$. By \cite[Proposition 1.10]{Hum}, we have $l(s_iw) =
l(vu) = l(v) + l(u)$.  But we also have $wu^{-1} = s_iv$ so that
$l(s_iv) = l(w) + l(u)$.  Suppose first that $l(s_iw) = l-1$. Then
$l(v) = l-1 -l(u)$ and $l + l(u) = l(s_iv) \leq l - l(u)$ which
implies that $u = 1$ so $s_iw \in W^I$.  Now suppose that $l(s_iw)
= l+1$. Then we have $l - l(u) \leq l + l(u) \leq l + 2 - l(u)$
with equality holding for exactly one inequality.  If $l - l(u) =
l + l(u)$ then again we have $u = 1$.  Otherwise, $l(u) = 1$.

In the last case we have $l(v) = l(w)$ and $l(s_iw) = l(s_iv) =
l+1$.  Let $u = s_r$ for some simple generator $s_r$. By the
Strong Exchange Condition (\cite[Theorem 5.8]{Hum}), $v =
(s_iw).s_r$ is obtained from $s_iw$ by taking a reduced word
$s_is_{a_1}\cdots s_{a_l}$ of $s_iw$ and omitting one generator.
If that simple generator is the first $s_i$ then $v = w$ and we
are done.  Otherwise $v = s_i s_{a_1}\cdots \hat{s}_{a_j} \cdots
s_{a_l}$ where $\hat{s}_{a_j}$ denotes omission.  But then it is
clear that $l(s_iv) = l-1$, a contradiction.
\end{proof}


\section{Skew affine Schur functions}
\label{sec:skew} The action of $\tS_n$ on the set of $n$-cores
induces another representation of $\U_n$.  Let $\U_n$ act on
$\C[\p^n]$ by
\[
u_i \cdot \nu = \begin{cases} s_i\cdot \nu & \mbox{if $s_i\cdot
\nu$ is obtained from $\nu$ by adding boxes.} \\ 0 &
\mbox{otherwise}.
\end{cases}
\]
The fact that this defines an action of $\U_n$ is easy to verify.
In fact we have

\begin{prop}
\label{prop:tworep} The above action of $\U_n$ on $\C[\p^n]$ is
isomorphic to the action of $\U_n$ on $\C[\Snn]$ where for $w \in
\Snn$ we define
\[
u_i \cdot w = \begin{cases} s_i \cdot w & \mbox{if $s_iw \in \Snn$
and $l(s_iw) > l(w)$}  \\
0 & \mbox{otherwise}.
\end{cases}
\]
Thus the action of $\U_n$ on $\C[\Snn]$ is obtained from the
action on $\C[\tS_n]$ by setting to 0 all elements $w \notin
\Snn$. The isomorphism is given by identifying $\ll \in \p^n$ and
$w(\ll) \in \p^n$.
\end{prop}

More generally, one can define an action of $\U_n$ on $\C[S^J]$
for other parabolic subgroups $S_J$ of $\tS_n$.

\begin{proof}
It is straightforward to check that the formulae of the
proposition do define a representation of $\U_n$ on $\C[\Snn]$.
By Proposition \ref{prop:ncores}, for $\nu \in \p^n$ the $n$-core
$s_i\cdot \nu$ is always obtained from $\nu$ by either adding
boxes or removing boxes or doing nothing. Let $w = w(\nu)$. Then
by Lemma \ref{lem:cox} applied to $W = \tS_n$ and $W^I = \Snn$, we
have either $s_iw = w(\mu)$ for some $\mu = s_i\cdot \nu$ or $s_iw
\in w\Sn$. In the latter case, $s_i \cdot \nu = \nu$.  In the
former case, using Proposition \ref{prop:Las}, adding boxes
corresponds to the case that $l(s_iw)
> l(w)$.
\end{proof}

Equip $\C[\p^n]$ with the inner product $\ip{\nu,\mu} =
\delta_{\nu\mu}$.

\begin{definition}
Let $\mu \subset \nu$ be two $n$-cores such that there is some $w
\in \tS_n$ satisfying $u_w \cdot \mu = \nu$.  The \emph{skew
affine Schur function} $\tF_{\nu/\mu}(X)$ is given by
\[
\tF_{\nu/\mu}(X) = \sum_{a = (a_1,a_2,\ldots,a_t)}
\ip{h_{a_t}(u)h_{a_{t-1}}(u) \cdots h_{a_1}(u) \cdot \mu, \nu}
x_1^{a_1} x_2^{a_2} \cdots x_t^{a_t}.
\]
\end{definition}

Suppose $\mu, \nu \in \p^n$.  Then using Proposition
\ref{prop:tworep} there is at most one permutation $w$ satisfying
$u_w \cdot \mu = \nu$. If $v$ and $u$ are right-Grassmannian
permutations corresponding to $\mu$ and $\nu$ then $w$ is given by
$w = uv^{-1}$ assuming that $l(w) + l(v) = l(u)$.  By Proposition
\ref{prop:Vaffine}, we have $\tF_{\nu/\mu} =\tF_w$ so that skew
affine Schur functions are special cases of affine Stanley
symmetric functions.  We write $w = w(\nu/\mu)$.  It is not true
that $l(w)$ is equal to the number of boxes in $\nu/\mu$, since
the action of $u_i$ may add more than one box.  It is also not
true that some $w$ exists satisfying $u_w \cdot \mu = \nu$ for
every pair of $n$-cores $\mu \subset \nu$ where containment is as
subsets of the plane.  For example, $(2,1,1) \subset (5,3,1)$ and
both are $3$-cores but such a $w$ does not exist.

When $\mu = \emptyset$, the permutation $w(\nu) =
w(\nu/\emptyset)$ is right-Grassmannian as defined earlier.  In
this case, the skew affine Schur function $\tF_{\nu/\emptyset}$ is
an affine Schur function. We write $\psi: \pa^n \rightarrow \p^n$
for the bijection satisfying $\tF_{\ll} =
\tF_{\psi(\ll)/\emptyset}$.

%

\section{Cores and $k$-tableaux}
\label{sec:ktableaux} One can view the skew affine Schur function
$\tF_{\nu/\mu}$ as the generating function for certain
semistandard tableaux built on $n$-cores.  These tableaux are
called \emph{$k$-tableaux} (with $k = n-1$) by Lapointe and Morse
\cite{LM04}. A (semistandard) $k$-tableau of shape $\nu/\mu$ and
weight $w(T) = (a_1,a_2,a_3,\ldots,a_l)$ is a chain of partitions
$\mu = \nu^{(0)} \subset \nu^{(1)} \subset \cdots \subset
\nu^{(l)} = \nu$ such that
\begin{enumerate}
\item
Each partition $\nu^{(i)}$ is an $n$-core.
\item
The successive differences $\nu^{(i)}/\nu^{(i-1)}$ contain at most
one box in each column.  That is they are horizontal strips.
\item
The contents $c(i,j) = j-i$ of the boxes $(i,j) \in
\nu^{(i)}/\nu^{(i-1)}$ involve exactly $a_i$ different residues
$\{r_1,r_2,\ldots,r_{a_i}\}$ modulo $n$ and $\nu^{(i)}$ has no
addable-corners with content congruent to one of these residues
$r_j$.
\end{enumerate}

\begin{figure}
\pspicture(0,20)(100,80) \abox(1,3){1} \abox(0,2){1} \abox(2,3){1}
\abox(3,3){1} \abox(4,3){2} \abox(5,3){2} \abox(6,3){3}
\abox(7,3){3} \abox(1,2){2} \abox(2,2){2} \abox(3,2){3}
\abox(4,2){3} \abox(0,1){3} \abox(1,1){3}

\endpspicture
\label{fig:k-tableau} \caption{A $k$-tableau with shape
$(5,4,2)/(1)$ and weight $(3,2,2)$.  Here $n = 4$.}
\end{figure}

When a $k$-tableau is drawn, the boxes of $\nu^{(i)}/\nu^{(i-1)}$
are filled with the number $i$.    We have (see also \cite{LM04})
\begin{prop}
\label{prop:ktableau} Let $\mu \subset \nu$ be two $n$-cores such
that there is some $w \in \tS_n$ satisfying $u_w \cdot \mu = \nu$.
Then
\[
\tF_{\nu/\mu}(X) = \sum_T x^{w(T)}
\]
where the sum is over all $k$-tableaux of shape $\nu/\mu$.
\end{prop}
\begin{proof}
The $n$-cores $\nu^{(i)}$ of a $k$-tableau are obtained by
successive applications of terms of $h_{a_i}(\u)$.  Thus
$\nu^{(i)} = u_{A_i} \cdot \nu^{(i-1)}$ for some term $u_{A_i}$ in
$h_{a_i}(\u)$.  This is equivalent to the description of
$k$-tableaux given above.  Condition (2) in the definition comes
from the fact that $u_{i+1}$ always precedes $u_i$ in the
definition of $h_k(\u)$ so that a box on a diagonal congruent to
$i$ modulo $n$ is never added after a box on a diagonal congruent
to $i+1$ modulo $n$.  Condition (3) follows from the description
(in Proposition \ref{prop:ncores}) of the action of $u_i$ on a
$n$-core, which adds all possible boxes along diagonals with
residue $i$.  The set $A_i$ is exactly the set of residues
$\{r_1,r_2,\ldots,r_{a_i}\}$.
\end{proof}

If $\ll$ is a partition fitting inside a $m \times (n-m)$ box for
some $m$ then at most $n-1$ diagonals are involved in $\ll$ and
necessarily $\ll \in \p^n$.  Within the $m \times (n-m)$ box, only
at most one box is added by the action of $s_i$.  In this case the
definition of a $k$-tableau reduces to a usual semistandard Young
tableau.  The following is then immediate.
\begin{prop}
\label{prop:schur} Let $\ll \subseteq ((n-m)^{m})$ for some $1
\leq m \leq n-1$. Then $\ll \in \pa^n \cap \p^n$ and
$\tF_{\ll/\emptyset} = \tF_{\ll} = s_\ll$.
\end{prop}

\section{Affine Schur and $k$-Schur functions}
\label{sec:kschur}

We now describe the relationship between affine Schur functions
and the $k$-Schur functions $\{s^{(k)}_\ll(X;t)\}$ (with $k =
n-1$).  The $k$-Schur functions $\{s^{(k)}_\ll(X;t)\}$ form a
basis of $\nL(t)$ and were originally used to investigate
Macdonald polynomial positivity.  Let $H_\mu(X;q,t)$ be given by
the plethystic substitution $H_\mu(X;q,t)= J_\mu(X/(1-q);q,t)$
where $J_\mu(X;q,t)$ is the integral form of Macdonald polynomials
\cite{Mac}.  Let $K^{(k)}_{\nu\mu}(q,t)$ and
$\pi^{(k)}_{\ll\nu}(t)$ be given by
\[
H_\mu(X;q,t) = \sum_\nu K^{(k)}_{\nu\mu}(q,t) s^{(k)}_\nu(X;t)
\;\;\; ; \;\;\; s^{(k)}_\nu(X;t) = \sum_\ll \pi^{(k)}_{\ll\nu}(t)
s_\ll(X).
\]
Then it is conjectured that $K^{(k)}_{\nu\mu}(q,t) \in \N[q,t]$
and $\pi^{(k)}_{\ll\nu}(t) \in \N[t]$ which would refine the
(proven) ``Macdonald positivity conjecture'' that the Schur
expansion of $H_\mu(X;q,t)$ has coefficients in $\N[q,t]$; see
\cite{Hai}.

There are a number of different definitions of $k$-Schur functions
\cite{LLM,LM03} which conjecturally agree.  The definition of the
$k$-Schur functions that we will use is from \cite{LM04} and is
(conjecturally) the $t=1$ specialisations of the original
definitions but are usually still called $k$-Schur functions.
Suppose $\tF_{\ll}(X) = \sum_\mu K^{(n)}_{\ll\mu} m_\mu$ where
$\ll \in \pa^n$ and the sum is over $\mu \in \pa^n$.  Then using
Proposition \ref{prop:ktableau} and the results of \cite{LM04},
the $k$-Schur functions $s^{(k)}_\ll(X) \in \nL$ are given by
requiring that
\[
h_\mu(X) = \sum_\ll K^{(n)}_{\ll\mu} s^{(k)}_\ll(X).
\]
This definition is called the \emph{$k$-Pieri} rule.

\begin{prop}
\label{prop:dual} Affine Schur functions and $k$-Schur functions
are dual bases of $\nL$ and $\Ln$, so that
$\ip{s^{(k)}_\mu,\tF_\nu} = \delta_{\mu\nu}$. \end{prop}
\begin{proof}
Write the affine Cauchy kernel
\begin{align*}
\Omega^{(n)}(X,Y) &= \sum_{\mu : \; \mu \in \pa^n} h_\mu(X)
m_\mu(Y) = \sum_{\mu : \;\mu \in \pa^n} \left( \sum_{\ll : \; \ll
\in \pa^n} K^{(n)}_{\ll\mu}
s^{(k)}_\ll(X) \right) m_\mu(Y) \\
& =  \sum_{\ll : \; \ll \in \pa^n} s^{(k)}_\ll(X) \left(\sum_{\mu
: \; \mu \in \pa^n} K^{(n)}_{\ll\mu} m_\mu(Y) \right)  = \sum_{\ll
: \; \ll \in \pa^n} s^{(k)}_\ll(X) \tF_\ll(Y),
\end{align*}
which is equivalent to duality.
\end{proof}

\section{Cylindric Schur functions}
\label{sec:cylindric}

In \cite{Pos}, Postnikov introduced and studied \emph{cylindric
Schur functions}, which he showed were symmetric functions; see
also closely related work of Gessel and Krattenthaler \cite{GK}.
Postnikov studied a special subset of the cylindric Schur
functions in finitely many variables which he called \emph{toric
Schur polynomials}.  He showed that the expansion coefficients of
toric Schur polynomials in the basis of Schur polynomials were
equal to 3-point genus 0 Gromov-Witten invariants
$C_{\lambda\mu\nu}^d$ of the Grasmannian $Gr_{m,n}$.  The Gromov
Witten invariant $C_{\lambda\mu\nu}^d$ counts the number of maps
$f: \mathbb{P}^1 \rightarrow Gr_{m,n}$ whose image has degree $d$
and meets generic translates of the Schubert varieties
$\Omega_\ll$, $\Omega_\mu$ and $\Omega_\nu$ at three marked points
$p_1,p_2,p_3 \in \mathbb{P}^1$.  In particular, these coefficients
are positive. They are the multiplicative constants of the (small)
quantum cohomology ring $QH^*(Gr_{m,n})$ of the Grassmannian.

In general cylindric Schur functions do not expand positively in
terms of Schur functions. See \cite{McN} for a detailed discussion
of this.

A \emph{cyclindric shape} $\ll$ is an infinite lattice path in
$\Z^2$, consisting only of moves upwards and to the right,
invariant under the translation by a vector $(n-m,-m)$ for some $m
\in [1,n-1]$.  We denote the set of such cylindric shapes by
$\c^{n,m}$.  If $\ll, \mu \in \c^{n,m}$ are cylindric shapes so
that $\mu$ always lies weakly to the left of $\ll$, then $\ll/\mu$
is a \emph{cylindric skew shape}.  We write $\mu \subset \ll$.

\begin{definition}
A \emph{cylindric semi-standard tableau} of shape $\ll/\mu$ and
weight $a = (a_1,a_2,\ldots,a_l)$ is a chain $\mu = \ll^{(0)}
\subset \ll^{(1)} \subset \cdots \subset \ll^{(l)} = \ll$ of
cylindric shapes in $\c^{n,k}$ such that each
$\ll^{(i)}/\ll^{(i-1)}$ is a cylindric skew shape with at most one
box in each column and $a_i$ boxes in any $n$ consecutive columns.
\end{definition}
When we draw a cylindric semi-standard tableau, we place the
number $i$ into the boxes of $\ll^{(i)}/\ll^{(i-1)}$.  The columns
will then be strictly increasing and the rows weakly increasing
(see Figure \ref{fig:cylindric}).

\begin{definition}
Let $\ll/\mu$ be a cylindric skew shape.  Then the cylindric Schur
function $s^c_{\ll/\mu}$ is given by
\[
s^c_{\ll/\mu}(X) = \sum_T x^T
\]
where the sum is over all cylindric tableau $T$ of shape
$\ll/\mu$.
\end{definition}

\begin{figure}
\pspicture(0,20)(100,150) \rput(110,136){$\cdots$}  \bbox(4,8){1}
\bbox(5,8){1} \bbox(4,7){2}
\bbox(5,7){3}\bbox(4,6){3}\bbox(5,6){8}\bbox(4,5){4} \abox(2,6){1}
\abox(3,6){1} \abox(2,5){2}
\abox(3,5){3}\abox(2,4){3}\abox(3,4){8}\abox(2,3){4} \bbox(0,4){1}
\bbox(1,4){1} \bbox(0,3){2}
\bbox(1,3){3}\bbox(0,2){3}\bbox(1,2){8}\bbox(0,1){4}
\rput(-10,24){$\cdots$}
\endpspicture
\caption{A cylindric semi-standard tableau with $n= 4$ and $m =
2$.} \label{fig:cylindric}
\end{figure}

One can alternatively define cylindric Schur functions in the same
way as skew affine Schur functions by letting $\U_n$ act on
infinite bit sequences $p = (\ldots,
p_{-2},p_{-1},p_0,p_1,p_2,\ldots)$ satisfying the periodicity
condition $p_i = p_{i+n}$.  It is clear that periodic bit
sequences are closed under the action of $\tS_n$ and in fact form
$n+1$ finite orbits depending on the value of $m =  p_1 + p_2 +
\cdots + p_n \in [0,n]$.

If $\ll \in \c^{n,m}$ is a cylindric shape then $s_i\cdot \ll$ is
the cylindric shape obtained from $\ll$ by either adding boxes at
all corners along diagonals congruent to $i$ mod $n$, or removing
such boxes, or doing nothing.  Define $u_i: \C[\c^{n,m}]
\rightarrow \C[\c^{n,m}]$ by
\[
u_i \cdot \ll = \begin{cases} s_i \cdot \ll & \mbox{if $s_i\cdot
\ll$ is obtained from $\ll$ by adding boxes.} \\
0 & \mbox{otherwise.} \end{cases}
\]
This defines a representation of $\U_n$ on $\C[\c^{n,m}]$, and
equipping $\C[\c^{n,m}]$ with the natural inner product one can
check directly using the definition of cylindric semistandard
tableaux that for $\mu \subset \ll \in \c^{n,m}$ the function
$\tF^c_{\ll/\mu}$ given by
\[
\tF^c_{\ll/\mu}(X) = \sum_{a = (a_1,a_2,\ldots,a_t)}
\ip{h_{a_t}(u)h_{a_{t-1}}(u) \cdots h_{a_1}(u) \cdot \mu, \ll}
x_1^{a_1} x_2^{a_2} \cdots x_t^{a_t}
\]
is equal to the cylindric Schur function $s^c_{\ll/\mu}(X)$.

\begin{lem}
Suppose $\ll$ and $\mu$ are cylindric shapes.  Then there is at
most one $w \in \tS_n$ satisfying $u_w \cdot \mu = \ll$.
\end{lem}
\begin{proof}
Suppose $v$ and $w$ satisfy $u_w \cdot \mu = \ll$ and $u_v \cdot
\mu = \ll$.  Let $u_w = u_{a_l}u_{a_{l-1}}\cdots u_{a_1}$ and $u_v
= u_{b_j}u_{b_{j-1}}\cdots u_{b_1}$.  We may assume that $u_{a_1}
\neq u_{b_1}$ for otherwise we can reduce to a smaller case by
letting $\mu:= u_{a_1}\cdot \mu$.  So let the rightmost occurrence
of $u_i = u_{a_1}$ in $u_{b_j}u_{b_{j-1}}\cdots u_{b_1}$ be
$u_{b_r}$.  The cylindric shape $\mu$ must have an addable corner
along the $i$-th diagonal so in particular none of $u_{b_{r-1}},
u_{b_{r-2}}, \ldots, u_{b_1}$ is equal to $u_{i+1}$ or $u_i$ and
we can move $u_{b_r}$ to the right most position to get another
reduced word for $u_v$, and then reduce to a smaller case.
\end{proof}

By Propostition \ref{prop:Vaffine}, cylindric Schur functions are
thus also special cases of affine Stanley symmetric functions.  In
fact more is true.

\begin{prop}
\label{prop:cylindricSkewAffine} Every cylindric Schur function
$\tF^c_{\ll/\mu}$ is a skew affine Schur function.
\end{prop}
\begin{proof}
Let $w \in \tS_n$ satisfy $u_w \cdot \mu = \ll$. We show first
that there are generalised $n$-cores $\nu,\rho$ such that $u_w
\cdot \nu = \rho$, which immediately implies $\tF^c_{\ll/\mu} =
\tF_{\rho/\nu}$ (the definition of $\tF_{\rho/\nu}$ for
generalised $n$-cores is the obvious one). Here, a generalised
$n$-core is a $n$-core with the diagonal labels possibly shifted:
so if $p = (\ldots,p_{-2},p_{-1},p_0,p_1,\ldots)$ is the edge
sequence of a $n$-core then the sequence
$q=(\ldots,q_{-2},q_{-1},q_0,q_1,\ldots)$ given by $q_i:= p_{i+k}$
defines a generalised $n$-core.  Equivalently, generalised
$n$-cores are in bijection with offset sequences $\set{d_i \mid i
\in Z}$ satisfying $d_{i-n} = d_i + 1$.

The edge sequence $p(\nu)$ is obtained from $p(\mu)$ by setting
$p_N = 0$ for $N \geq n \cdot (l(w)+1)$ and $p_N = 1$ for $N \leq
-(n \cdot (l(w) + 1)$.  Since it is clear that $\nu$ is a
generalised $n$-core, $\rho = u_w \cdot \nu$ is also a generalised
$n$-core as long as it is non-zero.

So the ``central'' part of $p(\nu)$ looks the same as $p(\mu)$ and
the action of $\U_n$ on the central part is identical.  An entry
of the bit sequence is moved no more than one step for each action
by a simple generator, so in total it is moved no further than
$l(w)$ from its initial position.  The alteration of $p(\nu)$ is
thus sufficiently far away from the centre that the altered bits
cannot affect whether a box is added at each step of the action of
the simple generators of $w$ on $\nu$.  For the action of some
$u_i$ to be non-zero we need only ensure that $s_i$ adds a box
somewhere to the shape.

Finally, if $\nu$ and $\rho$ are two generalised $n$-cores with
the same ``shift'' given by $d_1(\nu) + d_2(\nu) + \cdots +
d_{n}(\nu) = d_1(\rho) + d_2(\rho) + \cdots + d_{n}(\rho)$ then
one can shift again to find genuine $n$-cores $\nu^+$ and $\rho^+$
so that $\tF_{\rho/\nu} = \tF_{\rho^+/\nu^+}$.
\end{proof}


\section{321-avoiding permutations}
\label{sec:321}
\begin{definition}
An affine permutation $w \in \tS_n$ is \emph{321-avoiding} if no
reduced word for $w$ contains a subsequence of the form $i (i+1)
i$.
\end{definition}
When $w \in S_n$, this definition is the same as $w$ ``avoiding''
the pattern $321$, as shown in \cite{BJS}. We can extend this
naturally to the affine case.

\begin{prop}
An affine permutation $w \in \tS_n$ is 321-avoiding if and only if
there do not exist indices $x < y < z \in \Z$ such that $w(x) >
w(y) > w(z)$.
\end{prop}
\begin{proof}
Suppose first that some reduced word for $w$ contains a
subsequence of the form $i(i+1)i$, so that $w = vs_i s_{i+1} s_i
u$.  Recall that $ws_i \gtrdot w$ if and only if $w(i) < w(i+1)$.
Let $v' = vs_is_{i+1}s_i$.  Since the word is reduced, we must
have $v(i) < v(i+1) < v(i+2)$ and $v'(i) > v'(i+1) > v'(i+2)$. But
since multiplying by each simple generator in $u$ increases the
length of the permutation, the 3 integers $a = v'(i), b = v'(i+1)$
and $c = v'(i+2)$ will never be swapped past each other again.  So
there are indices $x < y < z$ such that $w(x) = a, w(y) = b$ and
$w(z) = c$.

Conversely, suppose $w$ has three indices $x < y < z$ so that
$w(x) > w(y) > w(z)$.  We may assume that there is no index $t$ in
the open interval $(x,y)$ such that $w(t) > w(y)$ for otherwise we
can replace $x$ by $t$.  Similarly, there is no $r$ in $(y,z)$ so
that $w(r) < w(y)$.  Now if $x < y - 1$, we multiply $w$ by $s_x$
on the right where $x$ is to be taken modulo $n$ as usual. Let $w'
= ws_x$.  Since $w(x) > w(x+1)$, we have $l(w') = l(w) -1$.  Also
note that if $w'(z) \neq w(z)$ then we have $w'(z-1) = w(z)$. This
is because $w(z) < w(x)$ and $z > x$ so it is not possible that $z
= x + kn$ for some $k \in \Z$.  Similarly, $w'(y)$ can only have
been moved to the left compared to $w(y)$, so that it is never
moved past $w(z)$.  So there are indices $x+1 = x'<y'<z'$ so that
$w'(x') > w'(y') > w'(z')$.  Furthermore $z'-x' < z - x$.
Repeating this (also with the roles of $z$ and $x$ swapped) we
eventually obtain $w'' \in \tS_n$ and $y'' \in \Z$ so that
$w''(y''-1) > w''(y'') > w''(y''+1)$. Clearly, $w''$ is not
321-avoiding and since at each step going from $w$ to $w''$ the
length is reduced, some reduced word for $w$ contains a reduced
word for $w''$ as a subword.  This shows that $w$ is not
321-avoiding.

\end{proof}

\begin{thm}
\label{thm:321} Let $w \in \tS_n$ be 321-avoiding.  Then $\tF_w$
is equal to a cylindric Schur function (and thus by Proposition
\ref{prop:cylindricSkewAffine} also a skew affine Schur function).
\end{thm}
\begin{proof}
We proceed by induction on $l = l(w)$, the case $l(w) = 1$ being
trivial. So assume $w = s_i\cdot v$ with $l(v) = l(w) - 1$ and
that $u_v \cdot \mu = \nu$ for cylindric shapes $\mu,\nu$.  Pick a
reduced word $\rho = \rho_1\rho_2\cdots \rho_{l-1}$ for $v$. Pick
$k$ minimal so that $\rho_k = i$, if such a $k$ exists. Then since
$w$ is 321-avoiding, we must have unique $x,y < k$ satisfying
$\rho_x = i+1$ and $\rho_y = i-1$. We claim that $s_i \cdot \nu$
is obtained from $\nu$ by adding boxes.  This is clear since after
applying $s_{i+1}$ and $s_{i-1}$, the shape $\nu$ must have edge
sequence satisfying $p_i(\nu) = 0$ and $p_{i+1}(\nu) = 1$.

If no such $k$ exists and $p_i(\nu) = 1$ or $p_{i+1}(\nu) = 0$
then in the first case $i-1$ does not occur in $\rho$ and
$p_i(\nu) = p_i(\mu)$.  In the second case $i+1$ does not occur in
$\rho$ and $p_{i+1}(\nu) = p_{i+1}(\mu)$.  In either or both
cases, we let $\ll$ be the cylindric skew shape obtained from
$\mu$ by setting $p_{i+kn}(\ll) = 0$ and $p_{i+1+kn}(\ll) = 1$
(and keeping the rest of the edge sequence the same).  Then it is
clear that $u_w \cdot \ll \neq 0$ so that $\tF_w = \tF^c_{(u_w
\cdot \ll) /\ll}$.
\end{proof}

If $\ll$ and $\mu$ are cylindric shapes satisfying $u_w \cdot \mu
= \ll$ then $w$ is necessarily 321-avoiding.  In fact, the action
of $\U_n$ on cylindric shapes always satisfies the additional
relation $u_i u_{i+1} u_i = u_{i+1} u_i u_{i+1} = 0$.  However,
this is not true for $n$-cores.  For example, let $n = 3$ and $\mu
= (1)$. Let $w = s_1 s_2 s_1 = s_2 s_1 s_2$.  Then $w \cdot \mu =
(3,1,1)$.  This shows that skew affine Schur functions are
considerably more complicated than cylindric Schur functions.  In
fact more is true:

\begin{prop}
\label{prop:Sta} There exists $\mu \in \p^n$ so that for each $w
\in \Sn$, there is a $n$-core $\ll$ so that $F_w = \tF_{\ll/\mu}$.
\end{prop}
\begin{proof}
We can pick $\mu$ to be any $n$-core with offsets satisfying
$d_1(\mu) < d_2(\mu) < \cdots < d_n(\mu)$.  Then by Proposition
\ref{prop:ncores}, $u_w \cdot \mu \neq 0$ so that $F_w = \tF_{w
\cdot \mu /\mu}$.
\end{proof}


\section{Positivity}
\label{sec:positivity} We conjecture that affine Schur functions
generalise Schur functions for Stanley symmetric function
positivity (Theorem \ref{thm:pos}).

\begin{conjecture}
\label{conj:pos} The affine Stanley symmetric functions $\tF_w(X)$
expand positively in terms of the affine Schur functions
$\tF_\ll(X)$.
\end{conjecture}

This conjecture seems to be consistent with all the known
behaviour of $k$-Schur functions and cylindric Schur functions.

It has been conjectured \cite{LLM, LM03} that the multiplicative
constants $d^\ll_{\nu\mu}$ for $k$-Schur functions given by
\[s^{(k)}_\nu s^{(k)}_\mu = \sum_{\ll \in \pa^n} d^\ll_{\nu\mu}
s^{(k)}_\ll\] are non-negative.  In \cite{LM05}, it is shown that
the coefficients $d^\ll_{\nu\mu}$ include the multiplicative
constants of the Verlinde algebra of $U(m)$ at level $n-m$.
\begin{prop}
\label{prop:kschurmult} Conjecture \ref{conj:pos} implies
$d^\ll_{\nu\mu} \geq 0$.
\end{prop}
\begin{proof}
Using Proposition \ref{prop:dual}, together with equation
(\ref{eq:hopf}), we have
\[d^\ll_{\nu,\mu} = \ip{s^{(k)}_\nu s^{(k)}_\mu,\tF_\ll} = \ip{s^{(k)}_\nu \otimes
s^{(k)}_\mu, \Delta \tF_\ll}.\] But $\Delta \tF_\ll = \sum_{\rho
\subset \ll} \tF_{\ll/\rho} \tF_\rho$ where the sum is over $\rho
\in \pa^n$ such $\psi(\rho) \subset \psi(\ll)$ as $n$-cores, and
such that $\tF_{\ll/\rho}:=\tF_{\psi(\ll)/\psi(\rho)}$ is defined
(the bijection $\psi: \pa^n \rightarrow \p^n$ was defined in
Section \ref{sec:skew}). Using Proposition \ref{prop:dual} again,
we have $d^\ll_{\nu,\mu} = \ip{s^{(k)}_\nu, \tF_{\ll/\mu}}$ which
would be positive if Conjecture \ref{conj:pos} is true.
\end{proof}

Call a cylindric skew shape $\ll/\mu$ where $\ll,\mu \in \c^{n,m}$
\emph{toric} if the \emph{toric Schur polynomial}
$\tF^c_{\ll/\mu}(x_1,x_2,\ldots,x_m)$ is non-zero \cite{Pos}. Then
Postnikov showed that the coefficients $C^\ll_{\nu/\mu}$ given by
\[
\tF^c_{\nu/\mu}(x_1,\ldots,x_m) = \sum_\ll
C^\ll_{\nu/\mu}s_\ll(x_1,\ldots,x_m)
\]
were Gromov-Witten invariants of $Gr_{m,n}$.  These coefficients
are known to be non-negative from their geometric definition, but
a combinatorial proof is still lacking.

\begin{prop}
Conjecture \ref{conj:pos} implies that $C^\ll_{\nu/\mu} \geq 0$.
\end{prop}
\begin{proof}
Postnikov showed that the only Schur polynomials
$s_\ll(x_1,\ldots,x_m)$ which appear in the Schur expansion of
$\tF^c_{\nu/\mu}(x_1,\ldots,x_m)$ satisfy $\ll \subset ((n-m)^m)$.
By Proposition \ref{prop:schur}, these must be exactly the affine
Schur functions which occur in the affine Schur expansion of
$\tF^c_{\nu/\mu}(x_1,\ldots,x_m)$.
\end{proof}
See also McNamara's work on cylindric Schur positivity \cite{McN}.

\begin{remark}
By Proposition \ref{prop:cylindricSkewAffine} and the proof of
Proposition \ref{prop:kschurmult}, the coefficients
$C^\ll_{\nu/\mu}$ are special cases of multiplication coefficients
for $k$-Schur functions.  It is known \cite{Wit} that the Verlinde
algebra of $U(m)$ at level $n-m$ agrees with quantum cohomology of
$Gr_{m,n}$ at $q=1$. Thus our work shows that on the one hand the
connection between toric Schur functions and quantum cohomology
and on the other hand the connection between $k$-Schur functions
and the Verlinde algebra are equivalent.
\end{remark}

Since $s_\ll^{(k)}(X) \in \nL$ we have an element
$s_\ll^{(k)}(\u)\in \U_n$ (as before $k = n-1$).  The following
proposition is inspired by the paper of Fomin-Greene \cite{FG}.

\begin{prop}
\label{prop:kschur} Let $c_{w\ll} \in \Z$ be given by
\[
s_\ll^{(k)}(\u) = \sum_{w \in \tS_n} c_{w\ll} u_w.
\]
Then $c_{w\ll} = a_{w\ll}$ where $a_{w\ll}$ is the coefficient of
$\tF_\ll$ in $\tF_w$.
\end{prop}

\begin{proof}
We compute using the (non-commutative) affine Cauchy kernel that
\begin{align*}
\tF_w(X) = \sum_{\ll \in \pa^n} \ip{h_\ll(\u) \cdot 1, w} m_\ll(X)
= \ip{\Omega^{(n)}(x,\u) \cdot 1, w}        =  \sum_{\ll \in
\pa^n} \ip{s^{(k)}_\ll(\u) \cdot 1, w} \tF_\ll(X).
\end{align*}
Thus the coefficient of $\tF_\ll$ in $\tF_w$ is equal to
$c_{w\ll}$.
\end{proof}
Thus Conjecture \ref{conj:pos} is equivalent to $c_{w\ll} \geq 0$:
every non-commutative $k$-Schur function can be expressed as a
non-negative sum of monomials in $\{u_0,u_1,\ldots,u_{n-1}\}$.
When $\ll$ is contained in some $(n-m) \times m$ box, then the
$k$-Schur function $s^{(k)}_\ll$ is actually the Schur function
$s_\ll$ \cite{LM03}. If in fact $|\ll| \leq n-1$, then by
restricting to proper subsets of the generators
$\{u_0,u_1,\ldots,u_{n-1}\}$ (like in Proposition
\ref{prop:eformula}) one can give a positive monomial formula for
$s_\ll(\u)$ in terms of reading words of tableaux using the
results of \cite{FG} on non-commutative Schur functions. This for
example gives combinatorial interpretations of some Gromov-Witten
invariants corresponding to very small shapes. However, it is
likely that such combinatorial interpretations are easily obtained
from existing results.

\section{Final comments}
\label{sec:final}
\subsection{Which affine Stanley symmetric functions are Schur,
skew Schur or cylindric?} In \cite{BJS}, the question of which
Stanley symmetric functions equalled a skew Schur function was
studied.  As Proposition \ref{prop:Sta} indicates, the
corresponding problem for affine Stanley symmetric functions may
well be more difficult.  We call an affine permutation $w$
\emph{affine vexillary} (respectively \emph{skew affine vexillary}
or \emph{cylindric vexillary}) if $\tF_w$ is equal to some affine
Schur function (respectively some skew affine Schur function or
cylindric Schur function).

\begin{problem} \hspace{20pt}
\label{prob:vex} Which affine permutations are affine vexillary,
skew affine vexillary and cylindric vexillary?
\end{problem}

For example, Theorem \ref{thm:321} shows that all 321-avoiding
permutations are cylindric vexillary.  It is not clear whether
$\mu(w) = \ll(w)$ implies that $w$ is vexillary, in the notation
of Section \ref{sec:affineSchur}.  The corresponding statement is
true for usual permutations and follows from part (2) of Theorem
\ref{thm:Sta}.

Cylindric Schur and affine skew Schur functions arise from
representations of $\U_n$ on different sets of infinite bit
sequences.  It would be interesting to find other sets of infinite
bit sequences which are closed under the action of $\tS_n$ and to
define actions of $\U_n$ on them.

\subsection{The affine flag variety, quantum cohomology and fusion ring}
The connections with $k$-Schur functions and with cylindric Schur
functions indicate that affine Stanley symmetric functions are
important objects.

Our results show directly that $k$-Schur functions and cylindric
Schur functions are related.  In some cases, this was already
known if we combine Postnikov's work on cylindric Schur functions
and Gromov-Witten invariants of the Grassmannian with Lapointe and
Morse's work showing that multiplication $k$-Schur functions
calculate the multiplication in the fusion ring.  Finally it is
known that the fusion ring agrees with the quantum cohomology
$QH^*(Gr_{m,n})$ of the Grassmannian at $q = 1$ (\cite{Wit}).
These connections suggest that there may be an interesting
$q$-analogue of our theory.  It is not clear whether the
$q$-analogue in quantum cohomology is related to the $t$-analogue
of the original $k$-Schur functions $s^{(k)}_\ll(X;t)$ arising
from Macdonald polynomial theory.

However, the most interesting direction to take seems to be the
connections with the affine flag variety (type $A$).  Shimozono
has conjectured that the multiplication of $k$-Schur functions
calculate the homology multiplication of the affine Grassmannian.
The dual conjecture is that affine Schur functions represent the
Schubert classes in the cohomology of the affine Grassmannian
\cite{MS}.  There is computational evidence in support of these
conjectures. Extending these conjectures from the affine
Grassmannian to the affine flag variety would involve defining
\emph{affine Schubert polynomials} which should in some sense be
``unstable'' versions of affine Stanley symmetric functions.

\subsection{A dual version of $\tF_w$}
We have shown that affine Schur functions $\tF_\ll$ are dual to
the $k$-Schur functions $s^{(k)}_\ll(X)$.  The $k$-Schur functions
are conjectured to be Schur positive \cite{LM03} (in \cite{LLM}
the Schur positivity is part of the definition). Define the
\emph{dual affine Stanley symmetric function} $\tF^d_w$ by
\[
\tF^d_w(X) = \sum_{\ll} a_{w\ll} s^{(k)}_\ll(X)
\]
where as before $a_{w\ll}$ is given by $\tF_w = \sum_\ll a_{w\ll}
\tF_\ll$.  If Conjecture \ref{conj:pos} is true as well as the
Schur positivity of $k$-Schur functions, then $\tF_w$ would be
Schur positive.  If so, is it the character of a natural $S_m$ or
$GL(N)$ module?

\subsection{Affine stable Grothendieck polynomials}
Whereas Schubert polynomials are representatives for the
cohomology of the flag variety, \emph{Grothendieck polynomials}
are representatives for the K-theory of the flag variety.  In the
same way that Stanley symmetric functions are stable Schubert
polynomials, one can define \emph{stable Grothendieck
polynomials}.  Our definition of affine Stanley symmetric
functions naturally generalises to a definition of affine stable
Grothendieck polynomials (see \cite{FG} or \cite{FK}).

Let $\tilde{\U}_n$ be the algebra obtained from $\U_n$ by
replacing the relation $u_i^2 = 0$ with $u_i^2 = u_i$.  Define
$\tilde{h}_k(\u) \in \tilde{\U}_n$ for $k \in [0,n-1]$ with the
same formula as for $h_k(\u)$.

\begin{definition}
Let $w \in \tS_n$.  The affine stable Grothendieck polynomial
$\tilde{G}_w(X) \in \Ln$ is

\[\tilde{G}_{w}(X) = \sum_{a = (a_1,a_2,\ldots,a_t)}
\ip{\h_{a_t}(u)\h_{a_{t-1}}(u) \cdots \h_{a_1}(u) \cdot 1, w}
x_1^{a_1} x_2^{a_2} \cdots x_t^{a_t},\] where the sum is over
compositions of $l(w)$ satisfying $a_i \in [0,n-1]$.

\end{definition}
The functions $\tilde{G}_w$ are not homogeneous.  The highest
degree part of $\tilde{G}_w$ is equal to $\tF_w$.

\begin{thm}
\label{thm:groth} The affine stable Grothendieck polynomial
$\tilde{G}_w(X)$ is a symmetric function.
\end{thm}

We will first need the following lemma.
\begin{lemma}
\label{lem:groth} Let $a,b \in [0,n-1]$ satisfy $a \neq b + 1$.
Then in $\tilde{\U}_n$ we have
\[
u_{b}u_{b-1}\cdots u_a u_{b}u_{b-1}\cdots u_a = u_{b-1}\cdots u_a
u_{b}u_{b-1}\cdots u_a = u_{b}u_{b-1}\cdots u_a u_{b}u_{b-1}\cdots
u_{a+1}.
\]
\end{lemma}
\begin{proof}
The result follows easily by induction, the base case being the
defining identity $u_b^2 = u_b$.
\end{proof}

\begin{proof}[Proof of Theorem \ref{thm:groth}]
We show that $\tilde{h}_k(\u)\tilde{h}_l(\u) =
\tilde{h}_l(\u)\tilde{h}_k(\u)$, as in Proposition
\ref{prop:commute}.  Our approach will be the same as in
Proposition \ref{prop:commute}, but since not just reduced words
are involved, the proof is slightly more difficult.  We indicate
the modifications of the proof of Proposition \ref{prop:commute}
which are needed -- the global structure of the proof is
completely identical, but the calculation within each critical
interval is more delicate.  The main difference is that an outer
interval $A_i$ may overlap with its right neighbour $B_k$. Let
$A^*$ and $B^*$ be an outer interval and its right neighbour as
before.  We may no longer assume that $\min(A^*) = \max(B^*) +1$,
but nevertheless we define $U = \phi(A^* \cup_i T_i,B^* \cup_i
S_i) \subset [b,a]$ with a small modification.  So we begin with
$U = [b,a]$ and a changing index $i$ set to $i := a$ to begin
with. The index $i$ decreases from $a$ to $b$ and at each step the
element $i$ may be removed from $U$ according to the rule:
\begin{enumerate}
\item If $i \in A^*$ then we remove it from $U$
unless $i \in S_k$ or $i \in (A^* \cap B^*) + 1$ for some $k \in
[1,s]$.
\item If $i \in B^*$ and $i \notin A^*$ then we remove it from $U$ unless $i \in
T_k + 1$ for some $k \in [1,s]$.
\item Otherwise we do not remove $i$ from $U$ and set $i:= i-1$.
Repeat.
\end{enumerate}
When $|U| = d$ we stop the algorithm.  The proof follows
essentially as in Proposition \ref{prop:commute} but in addition
we need the following types of manipulations in $\tilde{\U}_n$ for
$i >p >  k > j > m > l$ (cyclically):
\begin{multline}
\label{eq:cool} (u_i u_{i-1} \cdots u_{j+1} u_j) (u_k u_{k-1}
\cdots u_j u_{j-1} \cdots u_{l+1} u_l) = \\ (u_i u_{i-1} \cdots
u_m)(u_k u_{k-1} \cdots u_{j+1} u_m u_{m-1} \cdots u_l),
\end{multline}
and
\begin{multline}
\label{eq:cool2} (u_i u_{i-1} \cdots u_{j+1} u_j) (u_k u_{k-1}
\cdots u_j u_{j-1} \cdots u_{l+1} u_l) = \\ (u_i u_{i-1} \cdots
u_p u_{k-1}u_{k-2} \cdots u_j)(u_p u_{p-1} \cdots u_l),
\end{multline}
which follow from Lemma \ref{lem:groth}.  So for example we have
$u_4 u_3 u_2 u_3 u_2 u_1 u_0 = u_4 u_3 u_2 u_1 u_3 u_1 u_0$. One
checks that (\ref{eq:cool}) and (\ref{eq:cool2}) are exactly the
relations needed at the ``overlap'' between $A^*$ and $B^*$ and
show that the definition of $U = \phi(A^* \cup_i T_i,B^* \cup_i
S_i)$ induces the desired bijection.  Unlike in Proposition
\ref{prop:commute} we cannot perform our calculations within the
affine symmetric group since some of our words are not reduced.
However, the arguments required are nearly identical, as the next
example should show.
\end{proof}

\begin{example}
We illustrate the map $U = \phi(A^* \cup_i T_i,B^* \cup_i S_i)$.
Suppose $[b,a] = [2,20]$ and $A^* = [13,20]$, $B^* = [2,14]$. Let
$S_1 = [16,18]$ and $T_1 = [8,11]$ and $T_2 = \set{5}$ be the
inner intervals.  Then $d = 13$ and $U =
\set{2,3,4,5,6,9,10,11,12,14,16,17,18}$.  We can compute that
\begin{multline*}
u_{T_1}u_{T_2} u_{A^*}u_{B^*} u_{S_1}  =
u_{T_1}u_{T_2}u_{[13,20]}u_{14}u_{[2,12]}u_{S_1} = \\
u_{T_2}u_{[13,20]}u_{14} u_{[2,12]}u_{S_1}u_{(T_1 + 1)} =
u_{T_2}u_{[6,20]}u_{14} u_{[2,6]} u_{(T_1+1)} u_{S_1}
\end{multline*}
so that $U' = [6,20] \cup \set{5}$.  Finally one checks that $B^*
\cup_i S_i = \phi(U',U)$.
\end{example}

When $w$ is 321-avoiding, then we obtain \emph{cylindric stable
Grothendieck polynomials} which should be related to the quantum
$K$-theory of the Grassmannian.

\end{document}